\documentclass[journal,12pt,onecolumn,draftclsnofoot,]{IEEEtran}
\usepackage{url}
\usepackage{makeidx}
\usepackage{graphicx}
\usepackage{amsmath}
\usepackage{latexsym}
\usepackage{epsfig}
\usepackage{color}
\usepackage[ruled]{algorithm2e}
\SetKwInput{KwData}{Data}
\usepackage{enumerate}
\usepackage[colorlinks,urlcolor=blue,citecolor=blue,linkcolor=blue]{hyperref}
\newtheorem{theorem}{Theorem}

\usepackage[caption=false,font=normalsize,labelfont=sf,textfont=sf]{subfig}
\usepackage{multirow}
\usepackage[switch,columnwise]{lineno}
\setpagewiselinenumbers
\modulolinenumbers[1]
\usepackage{amsfonts}

\DeclareMathOperator*{\argmin}{arg\,min}

\newcommand{\rev}[1]{{\color{black}#1}}
\newcommand{\revi}[1]{{\color{black}#1}}
\newcommand{\revii}[1]{{\color{black}#1}}

\newtheorem{cor}[theorem]{Corollary}

\expandafter\def\expandafter\normalsize\expandafter{%
	\normalsize
	\setlength\abovedisplayskip{1.8pt}
	\setlength\belowdisplayskip{1.8pt}
	\setlength\abovedisplayshortskip{1.8pt}
	\setlength\belowdisplayshortskip{1.8pt}
}

\allowdisplaybreaks

\begin{document}
	
\title{A Stochastic Multi-Agent Optimization Framework for Interdependent Transportation and Power System Analyses}

\author{Zhaomiao~Guo,~\IEEEmembership{} Fatima~Afifah,~\IEEEmembership{}  Junjian~Qi,~\IEEEmembership{Senior Member,~IEEE}, Sina Baghali,~\IEEEmembership{Student Member,~IEEE} \vspace{-0.2in}
	\thanks{
		Z. Guo (guo@ucf.edu, corresponding author),   F. Afifah, and Sina Baghali are with the Department of Civil, Environmental and Construction Engineering, and Resilient, Intelligent, and Sustainable Energy Systems(RISES) cluster, University of Central Florida, Orlando, FL 32816, USA.}
    \thanks{
    	J. Qi is with the Department of Electrical and Computer Engineering, Stevens Institute of Technology, Hoboken, NJ 07030 USA.}
}
	

\maketitle
\begin{abstract}
	We study the interdependence between transportation and power systems considering decentralized renewable generators and electric vehicles (EVs). We formulate the problem in a stochastic multi-agent optimization framework considering the complex interactions between EV/conventional vehicle drivers, \revi{renewable}/conventional generators, and independent system operators, with locational electricity and charging prices endogenously determined by markets. We show that the multi-agent optimization problems can be reformulated as a single convex optimization problem and prove the existence and uniqueness of the equilibrium. To cope with the curse of dimensionality, we propose ADMM-based decomposition algorithm to facilitate parallel computing. Numerical insights are generated using standard test systems in transportation and power system literature.
	
	
\end{abstract}

\section*{Nomenclature}


{\noindent Sets and Indices}

\begin{itemize}[]
    \item $\boldsymbol{\Xi}$: set of uncertain capacity factors, indexed by $\boldsymbol{\xi}$
    \item $\mathcal{G(N,A)}/\mathcal{G(I,E)} $: transportation/power graphs, with $\mathcal{N, I}$ (indexed by $n, i$) being vertex sets and $\mathcal{A, E}$ (indexed by $a, e$) being edge sets
	\item $\mathcal{I}^S/\mathcal{I}^C/\mathcal{I}^T$: node sets of renewable/conventional generators/charging stations, indexed by $i$
    \item $\mathcal{I}_i$: node sets connecting to node $i$, indexed by $j$
	\item $\bar{\mathcal{R}}$/$\bar{\mathcal{S}}$ ($\mathcal {R/S}$): set of origins/destinations of conventional(electric) vehicles, indexed by $r, s$
\end{itemize}

{\noindent Parameters and Functions}
\begin{itemize}[]
    \item $\alpha$: coefficient for augmented Lagrangian function
    \item $\boldsymbol{\beta}/\epsilon$: coefficients/error terms in drivers' utility function
    \item $\nu$: iteration number
    \item $A$: node-link incidence matrix of transportation network
    \item $B$: investment budget for renewable generators
    \item $B_{ij}$: susceptance of transmission line $(i, j)$
	\item $C_i^{S, I}(\cdot)/C_i^{S,O}(\cdot)$: aggregated investment/operation cost functions of renewable generators at location $i$
	\item $C^C_i(\cdot)$: aggregated production cost functions of conventional generators at location $i$
	\item $E^{rs}$: O–D incidence vector of O–D pair $rs$ with +1 at origin and -1 at destination
	\item $e_{rs}$: average EV charging demand from $r$ to $s$
	\item $l_i^C, u_i^C$: lower and upper bound of conventional generation at location $i$ 
	\item $l_i$: power load (excluding charging load) at node $i$
	\item $P(\boldsymbol{\xi})$: probability measure of scenario $\boldsymbol{\xi}$
	\item $Q_{r}$: EVs initial \revi{quantity at location} $r$
	\item $\bar{q}_{rs}$: conventional vehicle travel demand from $r$ to $s$
	\item $tt_a(\cdot)$: link travel time function of link $a$
    \item $tt_{rs}$: equilibrium path travel time from $r$ to $s$
	\item $U_{rs}$: deterministic component of utility measures for drivers going from $r$ to $s$
	\item $v_a$: aggregated traffic flow on link $a$
	\item $z_i^{\nu}$: the renewable investment at $i$ at iteration $\nu$
\end{itemize}

{\noindent Variables}
\begin{itemize}[]
    \item $\gamma$: dual variables for non-anticipativity constraints
    \item $\theta_{i}$:  phase angle at node $i$
    \item $\lambda_i$: charging price at location $i$
    \item $\rho_i$: wholesale electricity price at location $i$
    \item $\tau_{n'}^{rs}$: dual variable of $rs$ flow conservation at node $n'$
	\item $d_i$: energy purchased by ISO at node $i$ 
	\item $f_{ij}/v_{a}$:  transmission/transportation link flow
	\item $g^S_i/g^C_i$:  electricity generation from \revi{renewable}/conventional generators at node $i$
	\item $p_i$: charging load at node $i$
	\item $q_{rs}$: EV travel demand from node $r$ to node $s$
	\item $u^S_i$:  capacity of \revi{renewable generator} at node $i$
	\item $x_a^{rs}$: link traffic flow on link $a$ that travels from node $r$ to node $s$
\end{itemize}



\section{Introduction}\label{sec:introduction}
\IEEEPARstart{T}{ransportation} and power systems are increasingly interdependent. The majority of modern transportation management, operations, and control systems are powered by electricity \cite{miles2016hurricane}; power system supplies, maintenance, and restoration (e.g. crew dispatch and mobile generators) rely on an efficient transportation network \cite{lei2016mobile}. The emergence of transportation electrification (e.g., electric vehicles (EVs) further strengthens the interdependence between these two systems \cite{wei2016robust}.

However, the majority of literature on planning and operation of transportation and power systems tackles these two systems separately, with outputs of one system served as exogenous inputs for the other \cite{wei2019interdependence}. Due to a large-body of literature from both transportation and power research communities, we only highlight representative ones for conciseness. For detailed review, one can refer to \cite{mwasilu2014electric, tan2016integration,rahman2016review, jang2018survey}. 

	
In power system literature, EVs with controllable charging and vehicle-to-grid (V2G) capabilities are considered as a part of a smart grid to improve power system efficiency and reliability \cite{tushar2015cost, liu2016ev, chung2018electric, liu2019optimal}. To estimate power demand, EVs' arrival and departure rates are typically assumed to be known as constants or probability distribution \cite{wei2019interdependence}. For example, \cite{melhorn2017autonomous} investigates the provision of reactive power from EVs to reduce the probability of power quality violation given EV plug-in and departure time; \cite{peng2017dispatching} includes a discussion on EVs participating in frequency regulation in the power system treating travel demand as given; \cite{wei2016robust} uses robust optimization and reformulates transportation equilibrium as charging demand uncertainty set for distribution network optimization; \cite{guo2019impacts} investigates the potential benefits of dynamic network reconfiguration to distribution network, with EVs' spatial-temporal availability and their charging demand estimated from transportation network models. 

In transportation literature, range anxiety of EV drivers and their charging behaviors have been integrated into transportation system modeling and charging infrastructure planning \cite{zheng2017traffic}, which typically treat power systems as exogenous. For example, \cite{lam2017coordinated} studies the coordinated parking problem of EVs to support V2G services, assuming that EV parking can be centrally coordinated and the EVs parking demand at each parking facility is given. The literature on charging infrastructure deployment in the transportation network \cite{guo2016infrastructure, chen2017deployment} typically treat locational electricity prices as exogenous variables influencing the facility choice decisions of EV drivers. 

Separating transportation and power system analyses is justifiable when EV penetration level is low and the feedback effects between transportation and power systems are negligible. With an increasing EV penetration level in the future, the EV travel and charging patterns will largely affect the spatial and temporal distribution of power demand as well as ancillary services availability; on the other hand,  charging costs, power availability, and ancillary service incentives from power system will influence the charging behaviors of EVs and impact transportation mobility. Therefore, alternative modeling strategies are urgently needed to capture the close couplings and feedback effects between the interdependent transportation and power systems. 

To address this issue, some recent studies aim to model transportation and power systems simultaneously using network modeling techniques \cite{wei2019interdependence}. Most of them have a bilevel structure where a central planner at the upper level optimizes system welfare;  traffic equilibrium and optimal power flow equations are presented at the lower level as constraints. For example, \cite{he2013optimal} determines the optimal allocation of public charging stations to maximize social welfare considering transportation and power transmission system equilibrium; \cite{he2016sustainability} aims to decide optimal charging prices to minimize system cost, including power loss in the distribution grid and travel time in the transportation network, with traffic equilibrium and distribution optimal power flow (OPF) as constraints; \cite{wei2017expansion} minimizes total power and transportation costs by choosing optimal power and transportation expansion decisions. The problem is linearized and reformulated as an exact mixed integer convex program\revi{; \cite{xie2020two} proposes a detailed model for transportation and power systems separately, where a third non-profit entity is responsible for the cooperation of the systems by minimizing the total social costs associated with both of them.}

However, transportation and power system planning are not centrally controlled by a single decision entity. Therefore, co-optimizing these two systems only provides a lower bound of the cost, which may not be achievable due to the decentralized nature of the decision making in power and transportation systems. In addition, mathematical programming with complementarity  constraints, whose mathematical properties and computation algorithms have been extensively investigated in the past \cite{dempe2002foundations, colson2005bilevel}, is still extremely challenging to solve, not to mention for large-scale transportation and power networks facing uncertainties. To model a decentralized transportation and power interaction and mitigate the computational challenges, \cite{amini2018optimal} uses a simulation-based approach to iteratively solve transportation at least cost vehicle routing problems and optimal power flow problems, and communicates locational marginal prices and charging demand between these two systems. However, the convergence of the proposed algorithm may not be guaranteed. This work is further extended in \cite{amini2019distributed} to consider power, transportation, and social systems as different layers coupled in an internet of thing (IoT) framework. In \cite{zhou2019dynamic}, the authors assume that each charging station has a renewable energy source installed at the same location, and the electricity price is determined based on forecasting the generation of these sources and the feed-in-tariff. However, \cite{amini2019distributed, zhou2019dynamic} ignore the power system operation, e.g., the power flow constraints.


\revii{Limited study considers endogenous flow-dependent travel time and charging prices for charging decision making in the context of decentralized multi-agent transportation and power system modeling. For example \cite{ huang2018economic, shukla2019multi, mao2020location} model charging location choice based on travel distance only. \cite{zhang2019yen} proposes a more realistic approach for modeling the path selection of EV drivers by allowing en-route charging. However, charging cost is not involved in the decision making. \cite{wei2017expansion, wang2018coordinated} provides a linearized transportation and power distribution system model for the EV charging station and transportation system expansion planning. The charging station selection of EVs is based only on the time of travel, without considering other factors such as the preferability of the charging station and charging price. \cite{he2013optimal,guo2019impacts} consider charging prices as well as flow-dependent travel time, but charging prices are exogenous.}

\revii{In summary, existing literature suffers from one or more of limitations, including centralized decision making, exogenous charging prices, deterministic agent decision making and simplified intra- and inter- systems interactions. Each of these limitations hampers the further investigation of interdependent transportation and power systems due to an increasing trend of transportation electrification, both in terms of planning and operation. First, only if EV behavior is properly modeled in an analytical framework, we can analyze how to plan and operate transportation infrastructure systems to influence the EV behavior and system interaction outcomes. Second, some of these characteristics are fundamental characteristics of both systems, such as uncertainty of renewable generation and decentralized decision making structure. Ignoring these characteristics will only provide a lower bound of system costs, which can cause significant bias for decision making and/or may not be implementable in reality. Third, one of the critical interdependence between transportation and power systems is charging demand. Modeling endogenous prices is able to capture a feedback loop of charging demand and charging costs in spatially dependent transportation and power networks. Modeling endogenous prices also allows for further incentive design of charging infrastructure investment and charging costs to optimize both transportation and power systems.}

\revii{Addressing the above mentioned limitations in a holistic modeling framework faces significant challenges from both modeling and computational perspectives. First, modeling decentralized decision making of power and transportation systems with endogenous prices typically leads to a highly non-convex problem, where equilibrium existence/uniqueness and algorithm convergence are generally not guaranteed. Second, the curse of dimensionality brought by high dimensional uncertain parameters, non-convex system interactions, and large-scale transportation and power systems requires coordination of novel modeling approaches and algorithm design.}

In this paper, we propose a multi-agent optimization framework to study the interactions between transportation and power systems considering uncertain renewable generation and EVs routing and charging location choices. \revii{More specifically, our main contributions are two-fold. First, we model decentralized key decision-makers with prices and travel time endogenously determined by the model in a unified framework, with equilibrium existence and uniqueness proved. Second, we propose an exact convex reformulation of the non-convex equilibrium problem based on strong duality that can lead to both scenario and system decomposition and parallel computing.}


The remainder of this paper is organized as follows. Section II introduces the mathematical formulation of each stakeholder's decision-making problem. Section III presents an algorithm based on convex reformulation and decomposition. Numerical experiments on small- and medium-scale test systems overlaying transportation and power networks are discussed in Section IV. Finally, we conclude in Section V with a summary of our results, contributions, and future extensions.
	
\section{Mathematical Modeling}\label{sec:models}
We explicitly model five types of stakeholders: \revi{renewable energy} investors, conventional generators, ISO, and EV/conventional vehicle drivers, who interact with each other in a unified framework. We assume a perfectly competitive market for the power supply and charging market, i.e. individual decision-makers in both supply and demand sides (e.g., generators, drivers) do not have market power to influence equilibrium prices through unilaterally altering his/her decisions. For example, power suppliers will decide their investment and/or production quantities to maximize their own profits, considering the locational electricity prices in the wholesale market. Electricity and charging prices are endogenously determined within the model. \revi{Since we assume perfectly competitive markets, our models are only applicable when each individual agent does not have significant market power to influence the market price.} 
\subsection{\revi{Renewable Energy} Investors Modeling}

In a perfectly competitive market, since prices are considered exogenous by investors, the investment decisions can be calculated by aggregating a large number of investors into a dummy investor, whose cost functions are aggregated costs for all investors \cite{guo2016infrastructure}. The decision-making problems of \revi{renewable energy} investors are formulated as two-stage stochastic programming, shown in (\ref{mod:sp}). 
    	\begin{subequations}\label{mod:sp}
		\begin{align}
		&\max_{\substack{\boldsymbol{u}^S, \boldsymbol{g}^S} \geq \boldsymbol{0}} &&  \mathbb{E}_{\boldsymbol{\xi}}\sum_{i \in \mathcal{I}^S}\left[ \rho_{i,\boldsymbol{\xi}} g_{i, \boldsymbol{\xi}}^{S} - C_i^{S,O}(g_{i, \boldsymbol{\xi}}^S)\right]- \sum_{i \in \mathcal{I}^S} C_i^{S,I}(u_i^S)
		\label{obj:sp_obj} \\
		&\text{\quad \ s.t.} && g_i^S(\boldsymbol{\xi}) \leq \xi_i u_i^S,\ \forall i \in  \mathcal{I}^S, \boldsymbol{\xi} \in \boldsymbol{\Xi} \label{cons:sp_capa}\\
		& && \sum_{i \in \mathcal{I}^S} c_i^S u_i^S \leq B \label{cons:sp_budg}
		\end{align}
	\end{subequations}

In the first stage, investors decide the \revi{renewable energy} generation capacity $\boldsymbol{u}^S$ considering the future uncertain \revi{renewable} capacity factors $\boldsymbol{\xi}$ at each location in the future, which is determined by uncertain factors such as \revi{renewable} radiance intensity, weather, and temperature \cite{xiao2014impact}. In the second stage, given a realization of \revi{renewable} generation uncertainties $\boldsymbol{\xi}$ and the first stage decision variable $\boldsymbol{u}^S$, \revi{renewable} generators determine the generation quantities $\boldsymbol{g}^S$. \rev{Our second-stage model is at hourly level.} We assume investors aim to maximize their long-term expected profits \cite{bruno2016risk}, which is calculated as total expected net revenue (revenue subtracting operational cost) in the second stage minus investment cost in the first stage, as shown in (\ref{obj:sp_obj}).  Without loss of generality, we assume $C_i^{S,I}(\cdot)$  has a quadratic form, with a positive quadratic coefficient to reflect an increasing marginal cost of land procurement and $C_i^{S,O}(\cdot)$ is a linear function \rev{\cite{guo2016infrastructure}}. Constraint (\ref{cons:sp_capa}) guarantees the power output of \revi{renewable generator} does not exceed the capacity $u_i^S$ times the uncertain \revi{renewable} intensity parameter $\xi_i$. Constraint (\ref{cons:sp_budg}) is the budget constraint for the total investment. \revi{Note that since we focus on long-term planning, model (\ref{mod:sp}) could be adapted to model different renewable energy resources, such as solar panels and wind turbines, whose total output capacity is uncertain and influenced by natural resource availability.} 
\subsection{Conventional Generators}

 For each scenario $\boldsymbol{\xi} \in \boldsymbol{\Xi}$, conventional generators solve the following optimization problem to determine their generation quantity $g_i^C, \forall i \in \mathcal{I}^C$.  Notice that the decision variables for conventional generators, ISO, and drivers are all scenario dependent, but we omit the notation $\boldsymbol{\xi}$ for brevity.
	\begin{subequations}\label{mod:cp}
		\begin{align}
		&\max_{\substack{\boldsymbol{g}^C} \geq \boldsymbol{0}} && \hspace{-1.5cm}\sum_{i \in \mathcal{I}^C} \big(\rho_i g_i^C - C^C_i(g_i^C) \big)\label{obj:cg_obj} \\
		&\text{\quad \ s.t.} &&\hspace{-1.5cm} g_i^C \leq u_i^C,\ \forall i \in  \mathcal{I}^C \label{cons:cg_capa}\\
		& && \hspace{-1.5cm}g_i^C \geq l_i^C,\ \forall i \in  \mathcal{I}^C \label{cons:cg_lower_bound}
		\end{align}
	\end{subequations}	
Objective (\ref{obj:cg_obj}) maximizes the profits of conventional generators at each scenario $\boldsymbol{\xi}$, which is calculated as total revenue $\sum_{i \in \mathcal{I}^C}\rho_i g_i^C$ subtracting total production costs $\sum_{i \in \mathcal{I}^C}C^C_i(g_i^C)$. We assume $C_i^{C}(\cdot)$ has a quadratic form, which is consistent with the settings in IEEE test systems\footnote{Source \url{https://matpower.org/docs/ref/matpower5.0/}}. Constraints (\ref{cons:cg_capa})--(\ref{cons:cg_lower_bound}) are the upper and lower bounds for power generation at each conventional generator location $i \in \mathcal{I}^C$. 

\subsection{ISO Modeling}

\rev{While conventional Level 1, Level 2, and DC fast charging infrastructure are connected to low voltage distribution network, recent proposal on 350kW--1MW ultra fast chargers 
may need to be connected to transmission or sub-transmission systems \cite{wei2019interdependence}. In this study we focus on inter-city travel with ultra-fast charging stations that are directly connected to the sub-transmission network. }

ISO monitors, controls, and coordinates the operation of electrical power systems. While ISO has many specific tasks, we focus on their daily operation to determine the power purchase and transmission plan to maximize system efficiency, which will implicitly determine locational marginal prices. 
Denote a power system as $\mathcal{G}^P = (\mathcal{I}, \mathcal{E})$. The ISO decision making can be described by model (\ref{mod:iso}). Objective function (\ref{obj:iso}) minimizes total energy purchasing cost from both renewable and conventional generators, $\sum_{i \in \mathcal{I}^S\cup\mathcal{I}^C}\rho_i d_i$,  minus total energy revenue from charging stations, $\sum_{i \in \mathcal{I}^T}\lambda_i p_i$.  Notice that when $p_i < 0$, ISO purchases $-p_i$ energy from charging station $i$, and objective function (\ref{obj:iso}) means minimizing total energy cost for power systems. Constraint (\ref{cons:iso_phas})--(\ref{cons:iso_capa}) are  power flow constraints, where (\ref{cons:iso_phas}) gives line flow patterns under DC power flow assumptions; (\ref{cons:iso_flow}) guarantees the power balance at each location $i$; and (\ref{cons:iso_capa}) is the transmission capacity constraint. Model (\ref{mod:iso}) guarantees that ISO will prefer purchasing power from cheaper generators and supplying power to more demanded charging locations, given grid topological and physical constraints. 
	\begin{subequations}\label{mod:iso}
		\begin{align}
		&\min_{\substack{\boldsymbol{p}, \boldsymbol{d}} \geq \boldsymbol{0}, \boldsymbol{\theta}, \boldsymbol{f}} && \sum_{i \in \mathcal{I}^S\cup\mathcal{I}^C}\rho_i d_i + \sum_{i \in \mathcal{I}^T}\lambda_i (-p_i)  \label{obj:iso} \\
		&\text{\quad \ s.t.} &&B_{ij}(\theta_i - \theta_j) = f_{ij},\ \forall (i,j) \in  \mathcal{E} \label{cons:iso_phas}\\
		& &&\sum_{j \in I_i}f_{ij} = d_i - p_i - l_i,\ \forall i \in  \mathcal{I} \label{cons:iso_flow}\\
		& &&-u_{ij}^{P} \leq f_{ij} \leq u_{ij}^{P},\ \forall (i,j) \in  \mathcal{E} \label{cons:iso_capa}
		\end{align}
	\end{subequations}
\subsection{Drivers Modeling}
	EV drivers need to determine their charging locations and travel routes in a transportation graph, denoted as $\mathcal{G(N, A)}$. EVs departing from $r \in \mathcal{R}$ select a charging station $s \in \mathcal{S}$, with the utility function defined in (\ref{eq:util}). A similar utility function has been adopted in previous literature \cite{guo2019impacts} to describe charging facility location choices, and can be extended to include other relevant factors based on evolving EV charging behaviors without interfering the fundamental modeling and computational strategies presented in this paper. Utility function (\ref{eq:util}) reflect the trade-off of EV drivers between four aspects: locational attractiveness $\beta_{0,s}$, travel time $-\beta_1 tt_{rs}$, charging cost $-\beta_2 e_{rs}{\lambda_s}$, and an error term $\epsilon$. \rev{For example, an EV driver may not choose a charging facility on the shortest path or with the lowest charging cost. Instead, he/she will balance locational preference, travel time, and charging costs.} \rev{Notice that we use average charging demand for each $rs$ pair. This modeling framework can be naturally extended to incorporate the variation of charging demands by creating dummy $rs$ pairs.} For home charging or workplace charging, where travelers have fixed destinations, the corresponding locational attractiveness (i.e. $\beta_{0,s}$) will dominate the other utility components. 
\begin{equation}
	    U_{rs} = \beta_{0,s} -\beta_1 tt_{rs} - \beta_2 e_{rs}\lambda_s + \epsilon \label{eq:util}
	\end{equation}
	Different assumptions on the probability distribution of $\epsilon$ result in different discrete choice models. In this study, we adopt a multinomial logit model, in which $\epsilon$ follows an extreme value distribution. The outputs of discrete choice models have two interpretations. First, discrete choice models calculate the probability of a  vehicle choosing from different destinations. Second, the results describe the traffic distribution to different destinations at an aggregated level.
	
	The utility function (\ref{eq:util}) \rev{partially} depends on travel time $tt_{rs}$, which is determined by the destination and route choices of all the drivers (including EVs and conventional vehicles). We assume conventional vehicles have known destination choices, with origin-destination flow $\bar{q}_{rs}$. The destination choice of EVs $q_{rs}$ and path travel time $tt_{rs}$ are coupled. On one hand, selection of charging location $s$ will increase the travel demands on certain paths from $r$ to $s$ and influence the travel time of the transportation network; on the other hand, path travel time will affect destination choices as travel time is a factor in the utility function (\ref{eq:util}). To capture these couplings, we adopt combined distribution and assignment (CDA) model \cite{sheffi1985urban} to model their destination choices and route choices, as shown in (\ref{mod:cda}). 
\begin{subequations}\label{mod:cda}
	\begin{align}
		& \min_{\boldsymbol{x},  \boldsymbol{q} \geq \boldsymbol{0} }
		& & & &\sum_{a \in \mathcal{A}} \int_{0}^{v_a} tt_a(u) \mathrm{d}u \nonumber\\
		& & & & & \hspace{-1.5cm}+ \frac{1}{\beta_1}\sum_{r \in \mathcal{R}}\sum_{s \in \mathcal{S}} q_{rs}\left(\ln q_{rs} - 1 + \beta_2e_{rs}\lambda_s - \beta_{0,s}\right) \label{obj:cda_obj}\\
		& \text{\quad \ s.t.} 
		& & & &  \boldsymbol{v} = \sum_{r \in \mathcal{R}, s \in \mathcal{S}}\boldsymbol{x}^{rs} + \sum_{r \in \bar{\mathcal{R}}, s \in \bar{\mathcal{S}}}\boldsymbol{\bar{x}}^{rs} \; \label{cons:cda_v_x}\\
		& & & (\boldsymbol{\tau}^{rs})& & A\boldsymbol{x}_{rs} = q_{rs}E^{rs}, \; \forall r \in \mathcal{R}, s \in \mathcal{S} \label{cons:cda_x_q}\\
        & & & & & A\boldsymbol{\bar{x}}_{rs} = \bar{q}_{rs}E^{rs}, \; \forall r \in \bar{\mathcal{R}}, s \in \bar{\mathcal{S}} \label{cons:cda_x_bar_q}\\
		& & & & &  \sum_{s \in \mathcal{S}} q_{rs} = Q_r, \forall r \in \mathcal{R}\label{cons:cda_q_d}
	\end{align}
\end{subequations}
	\revi{The Bureau of Public Roads (BPR) function \cite{us1964traffic} is used here to determine the time of travel $tt_a(\cdot)$}. $\sum_{a \in \mathcal{A}}\int_{0}^{v_a} tt_a(u) \mathrm{d}u$ in objective function (\ref{obj:cda_obj}) is the summation of the area under all the link travel cost functions $tt_a(\cdot)$, which is the total travel time $\sum_{a \in \mathcal{A}}v_a tt_a(v_a)$ subtracts externalities caused by route choices.   The second part consists of the entropy of traffic distribution $q_{rs}(\ln q_{rs} - 1)$ and utility terms in (\ref{eq:util}). Objective function (\ref{obj:cda_obj}) does not have a physical interpretation \cite{sheffi1985urban}, and it is a potential function constructed to guarantee the optimal solutions of (\ref{mod:cda}) are consistent with the first Wardrop principal (a.k.a. user equilibrium\footnote{The journey times in all routes used are equal and less than those that would be experienced by a single vehicle on any unused routes.}) \cite{wardrop1952some} and the multinomial logit  destination choice assumption.  These conditions can be guaranteed by sufficient and necessary Karush-Kuhn-Tucker conditions of model (\ref{mod:cda}) regarding $\boldsymbol{x}_{rs}$ and $\boldsymbol{q}$. For detail proofs, one can refer to \cite{sheffi1985urban}. Constraint (\ref{cons:cda_v_x}) calculates link flows by summing link flows of EVs and conventional vehicles over all origin and destination pairs. Constraints (\ref{cons:cda_x_q})--(\ref{cons:cda_x_bar_q}) are the vehicle flow conservation at each node for EV travel demand ${q}_{rs}$ and conventional vehicle travel demand $\bar{q}_{rs}$, respectively. Constraint (\ref{cons:cda_q_d}) guarantees the summation of EV traffic flow distribution to each $s$ equals to the total EV travel demand from $r$, $Q_r$.  The equilibrial travel time for each OD pair $rs$ can be calculated as $tt_{rs} \doteq \tau^{rs}_r - \tau^{rs}_s$, where $\tau^{rs}_i$ is the dual variable for constraint (\ref{cons:cda_x_q}).

Notice that an implicit assumption of model (\ref{mod:cda}) is that EV drivers will charge at their destinations. \revi{In \cite{guo2016critical}, an extension of combined distribution and assignment (CDA) model is introduced, denoted as generalized combined distribution and assignment model (GCDA), where vehicle can charge at their origins, destinations, or en-route. In addition, the GCDA model also consider both one-way and round-trip travel. The reason we do not include that model is to avoid over complication of the presented modeling framework and to keep a clear focus on developing a multi-agent power transportation interaction model considering renewable energy and effective computational strategies. Replacing CDA model with GCDA model is a relatively straightforward process and does not make a major difference to the modeling framework, and the proposed computational strategies in Section \ref{sec:alg} still apply.}


\subsection{Market Clearing Conditions}
	
The power purchased and supplied by ISO at each location $i$, $d_i$ and $p_i$, need to be balanced with locational power generation and charging demand, respectively, in a stable market. Otherwise, some of the supply or demand cannot be fulfilled and the market-clearing prices will be adapted accordingly. \rev{Since we aim to provide a steady-state modeling framework for long-term planning problems, real-time agent deviation from the equilibrium will not be considered and the hourly market clearing conditions }can be stated as (\ref{eq:equi}), where (\ref{eq:equi_gene}) guarantees that total energy purchased by ISO is equal to total energy generated at each location; \revi{equation (\ref{eq:equi_char}) enforces the balance between charging supply and charging demand of EVs. $i(s)$ denotes the node index in the power graph that charging location $s$ connects to.  Equation (\ref{eq:equi_char}) does not include regular power demand, denoted by $l_i$, which has been considered in the ISO power balancing constraint (\ref{cons:iso_flow})}. 

Locational prices of electricity $\rho_i$ and charging $\lambda_s$ can be interpreted as the dual variables for the market clearing conditions, respectively. \rev{Notice that we will have market clearing condition (\ref{eq:equi}) for each scenario $\boldsymbol{\xi}$, which could result in different prices for peak and off-peak hours.}
	\begin{subequations}
		\begin{align}
		& (\rho_i) && d_i = g_i^S + g_i^C, \; \forall i \in \mathcal{I}^S\cup\mathcal{I}^C \label{eq:equi_gene}\\
		& (\lambda_s) && \sum_{r \in \mathcal{R}}e_{rs}q_{rs} = p_{i(s)}, \; \forall s \in \mathcal{S} \label{eq:equi_char}
		\end{align}
		\label{eq:equi}
	\end{subequations}

\section{Computational Approach}
\label{sec:alg}
	The system formulations for each stakeholder need to be solved simultaneously, as the decision processes of one agent depend on the decisions of the others. For example, locational prices are determined by the collective actions of all agents, as described in the market clearing conditions in (\ref{eq:equi}). 
    
	
Since each stakeholder's optimization problem is convex with mild assumptions on the cost functions of travel, investment, and generation, one can reformulate each optimization problem as sufficient and necessary complementarity problems (CPs), and solve all the CPs together\cite{guo2017stochastic}. However, solving CPs directly is challenging because of the non-convexity and high dimensionality. \rev{The main computational contribution of this paper is to} propose an exact convex reformulation for this problem, as more formally stated in Theorem \ref{thm:refo_perf}.
	\begin{theorem}[Convex Reformulation]\label{thm:refo_perf}
	    If $tt_a(\cdot)$, $C_i^{S,I}(\cdot)$, $C_i^{S,O}(\cdot)$, and $C_i^{C}(\cdot)$  are convex, the equilibrium states of agents' interactions in perfectly competitive market, i.e., (\ref{mod:sp}), (\ref{mod:cp}), (\ref{mod:iso}), (\ref{mod:cda}), and market clearing conditions (\ref{eq:equi}) are equivalent to solving a convex optimization, as formulated in (\ref{mod:eq_perf}).
	\end{theorem}
	\begin{subequations}\label{mod:eq_perf}
		\begin{align}
        		& \underset{\substack{(\boldsymbol{u,g,p,}\\ \boldsymbol{d,x,q})\geq \boldsymbol{0}} }{\text{min}}
		& & \sum_{i \in \mathcal{I}^S} C_i^{S,I}(u_i^S) +  \mathbb{E}_{\boldsymbol{\xi}} \bigg[ \sum_{i \in \mathcal{I}^S}C_i^{S,O}(g_{i,\boldsymbol{\xi}}^S) \nonumber \\
		&&&+  \sum_{i \in \mathcal{I}^C} C^C_i(g_{i,\boldsymbol{\xi}}^C) + \frac{\beta_1}{\beta_2}\sum_{a \in \mathcal{A}} \int_{0}^{v_{a,\boldsymbol{\xi}}} tt_a(u) \mathrm{d}u \nonumber\\
		& & &+ \frac{1}{\beta_2}\sum_{r \in \mathcal{R}}\sum_{s \in \mathcal{S}} q_{rs,\boldsymbol{\xi}}\left(\ln q_{rs,\boldsymbol{\xi}} - 1 - \beta_{0,s}\right) 	\bigg]\label{obj:comb}\\
		& \text{\quad \ s.t.} 
		& &  (\ref{cons:sp_capa})-(\ref{cons:sp_budg}), (\ref{cons:cg_capa})-(\ref{cons:cg_lower_bound}), (\ref{cons:iso_phas})-(\ref{cons:iso_capa}), \nonumber\\
        &&&(\ref{cons:cda_v_x})-(\ref{cons:cda_q_d}), (\ref{eq:equi_gene})-(\ref{eq:equi_char}) \nonumber
		\end{align}
	\end{subequations}
    \IEEEproof See Appendix \ref{app:proofs}. \IEEEQED
    
    {\it Remark 1}: The intuition behind reformulation (\ref{mod:eq_perf}) is the reverse procedures of Lagrangian relaxation, where we move the penalty terms (e.g., $\rho_id_i$, $\rho_ig_i^S$, and $\rho_ig_i^C$) from the objective functions back to constraints (e.g. (\ref{eq:equi_gene}--\ref{eq:equi_char})). This convex reformulation allows us to apply alternating direction method of multipliers (ADMM) , which leads to decomposition with guaranteed convergence properties\cite{boyd2011distributed}, in contrast to heuristic diagonal methods. Furthemore, the existence and uniqueness of systems equilibrium is stated in Corollary \ref{cor:exis_uniq_syst_equi}.
    
    \begin{cor}[Existence and Uniquness of Systems Equilibrium]\label{cor:exis_uniq_syst_equi}
	    If $tt_a(\cdot)$, $C_i^{S,I}(\cdot)$, $C_i^{S,O}(\cdot)$, and $C_i^{C}(\cdot)$  are strictly convex functions, the system equilibrium exists and is unique. 
	\end{cor}
    \IEEEproof See Appendix \ref{app:proofs}. \IEEEQED
    
    To develop a decomposition-based algorithm, we first relax the investment decision $u_i^S$ to be scenario dependent and introduce non-anticipativity constraints\footnote{Decisions made before uncertainties being revealed should be not be measurable by a specific scenario $\boldsymbol{\xi}$.}, as shown in (\ref{eq:non-anti}). Then, problem (\ref{mod:eq_perf}) can be reformulated as (\ref{mod:eq_perf_non_anti}).
    \begin{equation}\label{eq:non-anti}
		(\gamma_i(\boldsymbol{\xi})) \quad u_{i,\boldsymbol{\xi}}^S = z_i, \; \forall i \in \mathcal{I}^S,  \boldsymbol{\xi} \in \boldsymbol{\Xi}
	\end{equation}

	\begin{subequations}\label{mod:eq_perf_non_anti}
		\begin{align}
		& \underset{\substack{(\boldsymbol{u,g,p,}\\ \boldsymbol{d,x,q})\geq \boldsymbol{0}} }{\text{min}}
		& & \mathbb{E}_{\boldsymbol{\xi}} \bigg[\sum_{i \in \mathcal{I}^S}(C_i^{S,I}(u_{i,\boldsymbol{\xi}}^S)+ C_i^{S,O}(g_{i,\boldsymbol{\xi}}^S)) \nonumber \\
		&&&+  \sum_{i \in \mathcal{I}^C} C^C_i(g_{i,\boldsymbol{\xi}}^C) + \frac{\beta_1}{\beta_2}\sum_{a \in \mathcal{A}} \int_{0}^{v_{a,\boldsymbol{\xi}}} tt_a(u) \mathrm{d}u \nonumber\\
		& & &+ \frac{1}{\beta_2}\sum_{r \in \mathcal{R}}\sum_{s \in \mathcal{S}} q_{rs,\boldsymbol{\xi}}\left(\ln q_{rs,\boldsymbol{\xi}} - 1 - \beta_{0,s}\right) 	\bigg]\label{obj:comb2}\\
		& \text{\quad \ s.t.} 
		& &  (\ref{cons:sp_capa})-(\ref{cons:sp_budg}), (\ref{cons:cg_capa})-(\ref{cons:cg_lower_bound}), (\ref{cons:iso_phas})-(\ref{cons:iso_capa}), \nonumber\\
        &&&(\ref{cons:cda_v_x})-(\ref{cons:cda_q_d}), (\ref{eq:equi_gene})-(\ref{eq:equi_char}), (\ref{eq:non-anti}) \nonumber
		\end{align}
	\end{subequations}
   Notice that (\ref{mod:eq_perf_non_anti}) can be decomposed by scenarios and by systems if we relax non-anticipativity constraints (\ref{eq:non-anti}) and charging market clearing conditions (\ref{eq:equi_char}), respectively. We propose a solution algorithm based on ADMM,  as summarized in Algorithm \ref{alg:ADMM}. \rev{Since Algorithm \ref{alg:ADMM} is an application of the general ADMM approach on a convex optimization problem, convergence is theoretically guaranteed \cite{boyd2011distributed}.} 
   
   We use an augmented Lagrangian approach to relax constraints (\ref{eq:equi_char}) and (\ref{eq:non-anti}), with dual variables $\boldsymbol{\lambda}$ and $\boldsymbol{\gamma}$, respectively. The decision variables of model (\ref{mod:eq_perf_non_anti}) can be divided into two groups, $(\boldsymbol{u}, \boldsymbol{g}, \boldsymbol{p}, \boldsymbol{d})$ and $(\boldsymbol{x}, \boldsymbol{q}, \boldsymbol{z})$, and be updated iteratively. When fixing $(\boldsymbol{x}, \boldsymbol{q}, \boldsymbol{z})$ and add augmented Lagrangian terms $\sum_{i \in \mathcal{I}^S}[\gamma_{i, \boldsymbol{\xi}}^{\nu} (u_{i, \boldsymbol{\xi}}^S - z_i^{\nu}) + \frac{\alpha}{2}||u_{i, \boldsymbol{\xi}}^S - z_i^{\nu}||_2^2]$ and $\sum_{i \in \mathcal{I}^T}[\lambda_{i, \boldsymbol{\xi}}^{\nu} (-p_{i, \boldsymbol{\xi}} + \sum_{r \in \mathcal{R}}q_{ri,\boldsymbol{\xi}}^{\nu})+ \frac{\alpha}{2}||-p_{i, \boldsymbol{\xi}} + \sum_{r \in \mathcal{R}}q_{ri,\boldsymbol{\xi}}^{\nu}||_2^2]$ in model (\ref{mod:eq_perf_non_anti}), we will have model (\ref{mod:eq_perf_non_anti_admm1}) to update $(\boldsymbol{u}, \boldsymbol{g}, \boldsymbol{p}, \boldsymbol{d})$; likewise, when $(\boldsymbol{u}, \boldsymbol{g}, \boldsymbol{p}, \boldsymbol{d})$ are fixed, we will have model (\ref{mod:eq_perf_non_anti_admm2}) to update $(\boldsymbol{x}, \boldsymbol{q})$, and the $\boldsymbol{z}$ updates can be derived analytically, see equation (\ref{eq:aver_u_z}), using the fact that $\sum_{\boldsymbol{\xi} \in \boldsymbol{\Xi}}\gamma_{i,\boldsymbol{\xi}}^{\nu} = 0$ \cite{fan2010solving}. Notice that the updates of both $(\boldsymbol{u}, \boldsymbol{g}, \boldsymbol{p}, \boldsymbol{d})$ and $(\boldsymbol{x}, \boldsymbol{q})$  can be decomposed by scenarios, as described in Step 1 and Step 2 in Algorithm \ref{alg:ADMM}. Step 3 updates the dual variables $\boldsymbol{\lambda}$ and $\boldsymbol{\gamma}$, with step size equal to augmented Lagrangian parameter $\alpha$, so that the dual feasibility with respect to Step 2 is guaranteed in each iteration (see \cite{boyd2011distributed}). Decomposing transportation and power systems also allows for taking advantage of sophisticated algorithms in respective domains, which is not the focus of this paper and will be left for future research.
   
	\resizebox{1\columnwidth}{!}{%
	\begin{minipage}{\columnwidth}
\begin{subequations}\label{mod:eq_perf_non_anti_admm1}
	\begin{align}
		& \underset{\substack{(\boldsymbol{u,g,p,}\\ \boldsymbol{d})\geq \boldsymbol{0}} }{\text{min}}
		& & \sum_{i \in \mathcal{I}^S}\bigg[C_i^{S,I}(u_{i, \boldsymbol{\xi}}^S)+ C_i^{S,O}(g_{i, \boldsymbol{\xi}}^S)\bigg] +  \sum_{i \in \mathcal{I}^C} C^C_i(g_{i, \boldsymbol{\xi}}^C) \nonumber \\
		&&& \hspace{-1cm}+ \sum_{i \in \mathcal{I}^S}\bigg[\gamma_{i, \boldsymbol{\xi}}^{\nu} (u_{i, \boldsymbol{\xi}}^S - z_i^{\nu}) + \frac{\alpha}{2}||u_{i, \boldsymbol{\xi}}^S - z_i^{\nu}||_2^2\bigg] \nonumber \\
		&&&\hspace{-1cm} + \sum_{i \in \mathcal{I}^T}\bigg[\lambda_{i, \boldsymbol{\xi}}^{\nu} (-p_{i, \boldsymbol{\xi}} + \sum_{r \in \mathcal{R}}q_{ri,\boldsymbol{\xi}}^{\nu})+ \frac{\alpha}{2}||-p_{i, \boldsymbol{\xi}} + \sum_{r \in \mathcal{R}}q_{ri,\boldsymbol{\xi}}^{\nu}||_2^2\bigg] \label{obj:comb_admm1}\\
		& \text{\quad \ s.t.} 
		& &  (\ref{cons:sp_capa})-(\ref{cons:sp_budg}), (\ref{cons:cg_capa})-(\ref{cons:cg_lower_bound}), (\ref{cons:iso_phas})-(\ref{cons:iso_capa}), (\ref{eq:equi_gene}) \nonumber
	\end{align}
\end{subequations}
\end{minipage}
}

\resizebox{\columnwidth}{!}{%
\begin{minipage}{\columnwidth}
\begin{subequations}\label{mod:eq_perf_non_anti_admm2}
	\begin{align}
		& \underset{(\boldsymbol{x, q} )\geq \boldsymbol{0} }{\text{min}}
		& & \frac{\beta_1}{\beta_2}\sum_{a \in \mathcal{A}} \int_{0}^{v_{a, \boldsymbol{\xi}}} tt_a(u) \mathrm{d}u + \frac{1}{\beta_2}\sum_{r \in \mathcal{R}}\sum_{s \in \mathcal{S}} q_{rs,\boldsymbol{\xi}}\left(\ln q_{rs,\boldsymbol{\xi}} - 1 - \beta_0^s\right)  \nonumber \\
		&&& + \sum_{i \in \mathcal{I}^T}\bigg[\lambda_{i, \boldsymbol{\xi}}^{\nu} (-p_{i, \boldsymbol{\xi}}^{\nu+1} + \sum_{r \in \mathcal{R}}q_{ri,\boldsymbol{\xi}})+ \frac{\alpha}{2}||-p_{i, \boldsymbol{\xi}}^{\nu+1} + \sum_{r \in \mathcal{R}}q_{ri,\boldsymbol{\xi}}||_2^2\bigg] 	\label{obj:comb_admm2}\\
		& \text{\quad \ s.t.} 
		& &  (\ref{cons:cda_v_x})-(\ref{cons:cda_q_d}) \nonumber
		\end{align}
	\end{subequations}
	\end{minipage}
}
\begin{algorithm}[ht]
	\KwResult{$\boldsymbol{u}^{\nu}, \boldsymbol{g}^{\nu}, \boldsymbol{p}^{\nu}, \boldsymbol{d}^{\nu}, \boldsymbol{x}^{\nu}, \boldsymbol{q}^{\nu}, \boldsymbol{z}^{\nu}, \boldsymbol{\lambda}^{\nu}, \boldsymbol{\gamma}^{\nu}$}
	initialization: $\boldsymbol{u}^{0}, \boldsymbol{g}^{0}, \boldsymbol{p}^{0}, \boldsymbol{d}^{0}, \boldsymbol{x}^{0}, \boldsymbol{q}^{0}, \boldsymbol{z}^{0}, \boldsymbol{\lambda}^{0}, \boldsymbol{\gamma}^{0}$, $\alpha=1$, $\varepsilon = 0.001$, $\nu = 1$, $\textrm{gap} = \infty$ \;
	define:  $\boldsymbol{y}_1=(\boldsymbol{u}^{\nu}_{\boldsymbol{\xi}}, \boldsymbol{g}^{\nu}_{\boldsymbol{\xi}}, \boldsymbol{p}^{\nu}_{\boldsymbol{\xi}}, \boldsymbol{d}^{\nu}_{\boldsymbol{\xi}})$,  $\boldsymbol{y}_2=(\boldsymbol{x}^{\nu}_{\boldsymbol{\xi}}, \boldsymbol{q}^{\nu}_{\boldsymbol{\xi}})$, $\mathcal{K}_1=\{\boldsymbol{y}_1|\boldsymbol{y}_1\ge 0, (\ref{cons:sp_capa})-(\ref{cons:sp_budg}), (\ref{cons:cg_capa})-(\ref{cons:cg_lower_bound}), (\ref{cons:iso_phas})-(\ref{cons:iso_capa}), (\ref{eq:equi_gene})\}$, $\mathcal{K}_2=\{\boldsymbol{y}_2|\boldsymbol{y}_2\ge 0, (\ref{cons:cda_v_x})-(\ref{cons:cda_q_d}) \}$ 

\While{$\textrm{gap} \geq \varepsilon$}{
		Step 1: $\forall \boldsymbol{\xi} \in \boldsymbol{\Xi}, \boldsymbol{u}^{\nu}_{\boldsymbol{\xi}}, \boldsymbol{g}^{\nu}_{\boldsymbol{\xi}}, \boldsymbol{p}^{\nu}_{\boldsymbol{\xi}}, \boldsymbol{d}^{\nu}_{\boldsymbol{\xi}} \in \underset{\boldsymbol{y}_1 \in \mathcal{K}_1}{\argmin} \; (\ref{obj:comb_admm1})$\;
	
		Step 2: $\forall \boldsymbol{\xi} \in \boldsymbol{\Xi}$, $\boldsymbol{x}^{\nu}_{\boldsymbol{\xi}}, \boldsymbol{q}^{\nu}_{\boldsymbol{\xi}} \in \underset{\boldsymbol{y}_2 \in \mathcal{K}_2}{\argmin}\; (\ref{obj:comb_admm2})$
		\begin{equation}\label{eq:aver_u_z}
		    \boldsymbol{z}^{\nu} = \sum_{\boldsymbol{\xi} \in \boldsymbol{\Xi}} P(\boldsymbol{\xi})\boldsymbol{u}^{\nu}_{\boldsymbol{\xi}}
		\end{equation}
	
	Step 3: 
	\begin{equation}
	    \forall i \in \mathcal{I}^S, \boldsymbol{\gamma}_{i, \boldsymbol{\xi}}^{\nu} = \boldsymbol{\gamma}_{i, \boldsymbol{\xi}}^{\nu-1} + \alpha (u_{i, \boldsymbol{\xi}}^{S,\nu} - z_i^{\nu})
	\end{equation}
	\begin{equation}
	\forall i \in \mathcal{I}^T, \boldsymbol{\lambda}_{i, \boldsymbol{\xi}}^{\nu} = \boldsymbol{\lambda}_{i, \boldsymbol{\xi}}^{\nu-1} + \alpha \bigg (-p_{i, \boldsymbol{\xi}}^{\nu} + \sum_{r \in \mathcal{R}}q_{ri,\boldsymbol{\xi}}^{\nu}\bigg)\end{equation}
	
	Step 4: 
	\vspace{-0.2\baselineskip}
	\begin{align*}
	    \textrm{gap}_1 &=  \max_{\boldsymbol{\xi} \in \boldsymbol{\Xi}, i \in \mathcal{I}^S} |u_{i, \boldsymbol{\xi}}^{S,\nu} - z_i^{\nu}|/z_i^{\nu}\\
	    \textrm{gap}_2 &=  \max_{\boldsymbol{\xi} \in \boldsymbol{\Xi}, i \in \mathcal{I}^T} |-p_{i, \boldsymbol{\xi}}^{\nu} + \sum_{r \in \mathcal{R}}q_{ri, \boldsymbol{\xi}}^{\nu}|/p_{i, \boldsymbol{\xi}}^{\nu}\\
	    \textrm{gap} &= \max\{\textrm{gap}_1, \textrm{gap}_2\}, \hspace{0.5cm} \text{let}: \nu := \nu + 1
	\end{align*}

	}
	\caption{ADMM-Based Decomposition Algorithm}\label{alg:ADMM}
\end{algorithm}

	\section{Numerical Simulation}
	
	\subsection{Three Node Test System}
	
\begin{figure}[!t]
    	\begin{minipage}[b]{0.5\linewidth}
            \centering
            \includegraphics[width=1\linewidth]{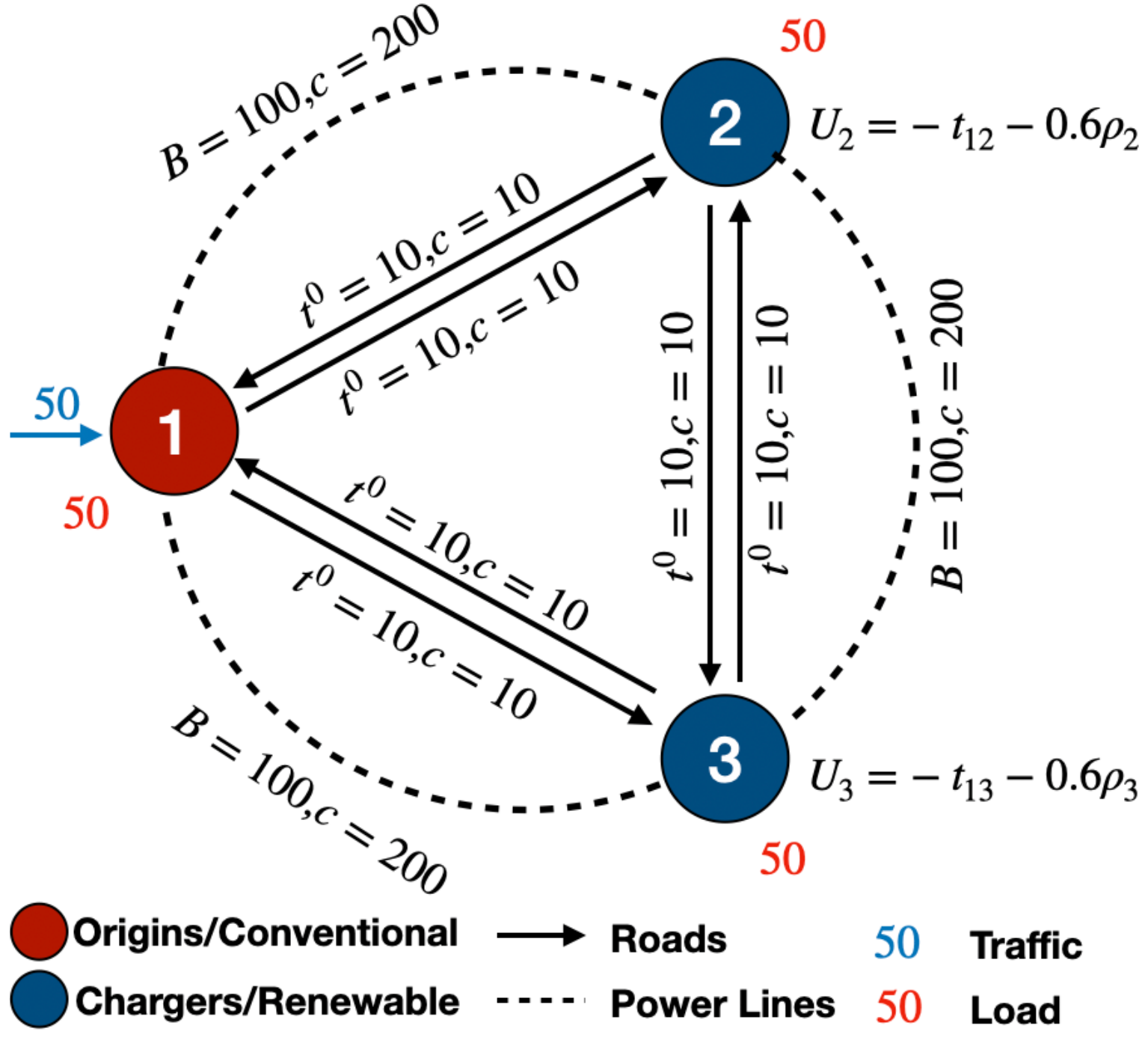}
            \caption{Three-node test system.}\label{fig:three_node}
        \end{minipage}
        \hspace{20pt}
        \begin{minipage}[b]{0.3\linewidth}
            \centering
            \includegraphics[width=1\linewidth]{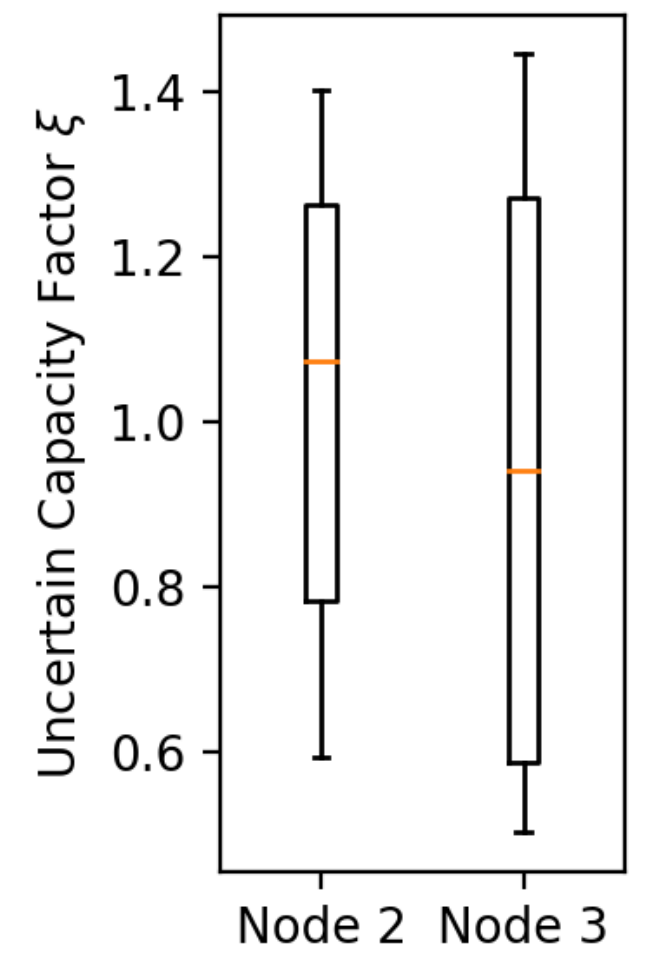}
            \caption{Variation in $\xi$.}
            \label{fig:sp_cap_scale}
        \end{minipage}
	\end{figure}

	To generate insights on the interdependence between transportation and power systems, we start with a three-node test system, as shown in Fig. \ref{fig:three_node}. The solid and dashed lines in Fig. \ref{fig:three_node} are transportation and power links, respectively. A traffic flow of 50 departing from node 1 has two possible charging destinations, nodes 2 and 3. Each node has 50 units of existing power load in addition to the charging load. All the other link parameters are shown in Fig. \ref{fig:three_node}. Fig. \ref{fig:sp_cap_scale} shows the distribution of uncertain renewable capacity factors for 10 scenarios. We implement our models and decomposition algorithms on Pyomo 5.6.6 \cite{hart2017pyomo} and solve the sub-problems using Cplex 12.8 and IPOPT 3.12.13, with 0.1\% optimality gap. All the numerical experiments presented in this paper were run on a 2.3 GHz 8-Core Intel Core i9 with 16 GB of RAM memory, under Mac OS X operating system\rev{, with parallel computing enabled}.

	We compare three cases. Case 1 is the base case, as shown in Fig. \ref{fig:three_node}. In case 2, the transportation capacity for link 1-3 reduces from 10 to 5. In case 3, the transmission capacity for link 1-3 reduces from 200 to 50. \revi{The reduction of transportation or transmission capacity may have multiple conflicting impacts on the system, including influence of equilibrium prices, renewable energy investment, system costs, and redistribution of traffic or power flow. Due to the interconnection and interdependence between transportation and power systems, these effects cannot be properly quantified without a modeling framework that can include transportation and power system decision making as endogenous variables and capture the feedback effects between these two systems. The results presented in this section aim to demonstrate the capability of our models on capturing the systems interdependence and quantifying the impacts of capacity reduction on system interaction outcomes.}

    \subsubsection{Effects on Equilibrium Prices}\label{sec:equi_price}
\begin{figure} 
        \begin{minipage}[b]{0.32\linewidth}
            \centering
            \subfloat[\footnotesize{Energy Prices}\label{fig:energy_prices}]{%
                \includegraphics[width=\linewidth]{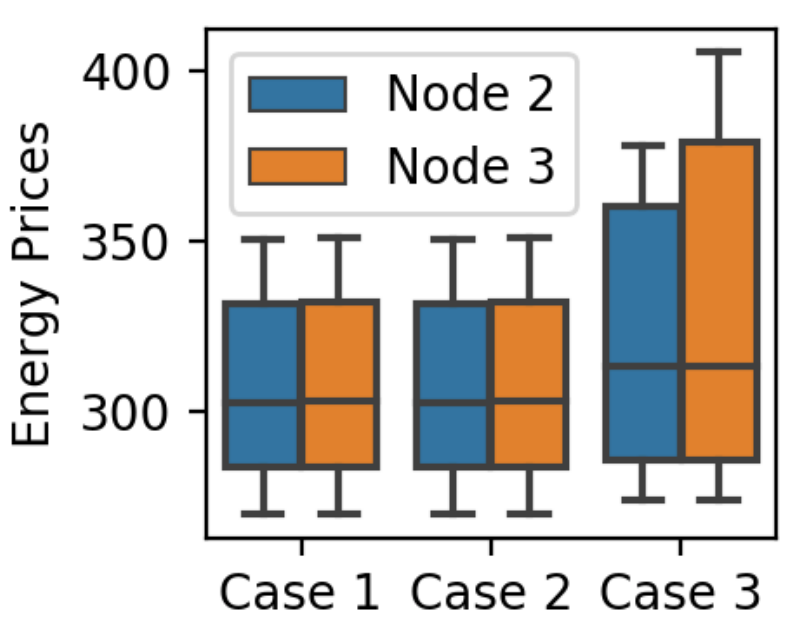}}
        \end{minipage}
         \begin{minipage}[b]{0.28\linewidth}
            \centering
            \subfloat[\footnotesize{Investment}\label{fig:investment}]{%
                \includegraphics[width=\linewidth]{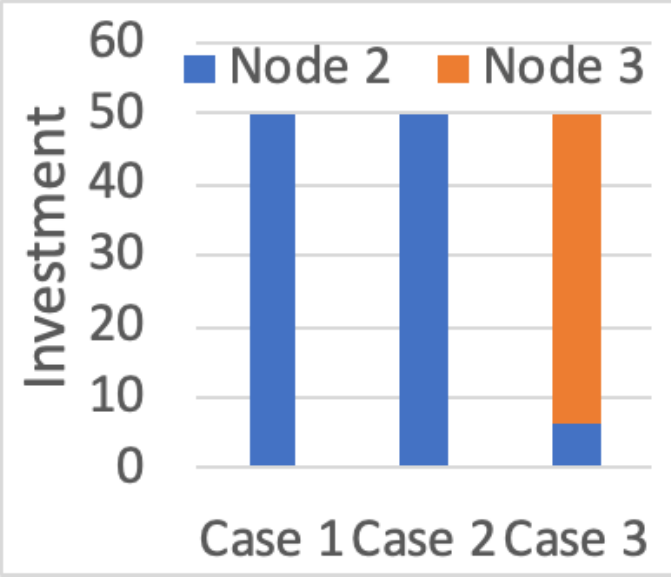}}
        \end{minipage}
        \begin{minipage}[b]{0.38\linewidth}
            \centering
            \subfloat[\footnotesize{Sytem Costs}\label{fig:system_cost}]{%
                \includegraphics[width=\linewidth]{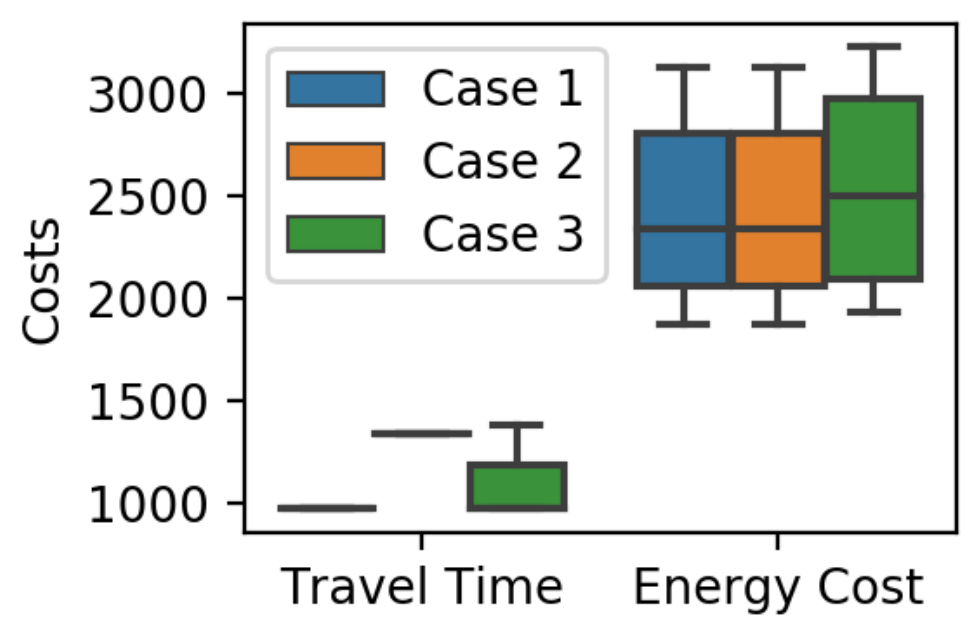}}
        \end{minipage}
        \caption{Impacts of link capacity on energy prices and system costs.}
        \label{fig:capacity_impact_1}
    \end{figure}

    From Fig. \ref{fig:energy_prices}, with sufficient transmission link capacity (cases 1 and 2), the energy prices at node 2  and 3 are identical, otherwise there will be incentives for the generators and ISO to supply more energy to the node with higher prices. \revi{Notice that the reduction of transportation link capacity does not impact the energy prices when there are sufficient transmission capacity and flexibility.} When we reduce transmission capacity for link 1-3, prices at both nodes 2 and 3 increase. The reason is that the transmission capacity of the whole network will not be better off with a reduction of transmission capacity on any links. \revi{With a limitation on overall transmission capacity, power may have to be generated at a more expensive locations, which leads to higher marginal prices.} But prices on node 3 increase to a larger extend than node 2, because node 3 is directly impacted by the transmission capacity limitation, which makes node 3 more challenging to receive energy. 

\subsubsection{Effects on Renewable Investment}
        
\revi{Since nodes 2 and 3 have the same cost parameters, the locational investment amount of renewable generators is determined by the equilibrium prices of energy $\boldsymbol{\rho}$ and the distribution of capacity factors $\boldsymbol{\xi}$, see model (\ref{mod:sp}). For example, a location with higher equilibrium price and a higher capacity factor will be strictly preferred. But in the case when these two factors are conflicting with each other, the investment amount will depend on whose influence is dominant.} The investment results are shown in Fig. \ref{fig:investment}. In cases 1 and 2, all investments are made in node 2. This is because node 2 has slightly higher capacity factors on average (see Fig. \ref{fig:sp_cap_scale}), and locational energy prices are the same at nodes 2 and 3 in cases 1, 2 (see Fig. \ref{fig:energy_prices}). \revi{But in case 3, the energy prices at node 3 is higher than the prices at node 2. The influence of energy prices on investment outperforms the influence of renewable capacity factors, which leads to more investment in node 3 in case 3. Comparing between cases 1 and 3, we can see that a reduction of transmission capacity for link  1-3 will increase the energy price at node 3 more than node 2, which will lead to a relocation of 46 units of investment \revii{from} node 2 to node 3. This observation indicates that renewable energy investment could offset some negative impacts of high energy costs and energy scarcity due to limited transmission capacity. Notice that we are only able to numerically quantify the influence of energy prices on renewable investment because prices are endogenously determined by the models and calculating the marginal impacts of prices on renewable investment based on model (\ref{mod:sp}) may be misleading.}

\subsubsection{Effects on System Costs}
System costs include travel costs and energy costs. Travel costs depend on traffic distribution and the capacity of transportation infrastructure. Fig. \ref{fig:system_cost} shows the distribution of travel costs and energy costs.  \revi{Both cases 2 and 3 have higher travel time than case 1, but for different reasons.} \revi{For case 2, the increased total travel time is because of a reduction on the transportation capacity on link 1-3, which leads to more congestion for the whole transportation network. The reason for higher travel time in case 3 than case 1 is because node 2 has a cheaper energy price compared to node 3 in case 3. More traffic will choose to travel through link 1-2. Since the total travel demand is fixed (i.e. 50), an imbalanced traffic pattern will cause more congestion.} Energy costs are the total costs of conventional energy and \revi{renewable} energy production. Since \revi{renewable} energy is cheaper and the total energy demand is fixed, cases with more \revi{renewable} energy utilized will have lower total energy costs. Cases 1 and 2 do not have transmission congestion, so all the \revi{renewable} energy can be utilized. The total energy in cases 1 and 2 will be identical and less than case 3, where transmission constraints prevent the effective usage of \revi{renewable} energy. \revi{These observations are consistent with common beliefs, which provides evidence that the modeling framework is effective to describe the main interaction between transportation and power systems.}
    
	\subsubsection{Effects on Flow Distribution}
\begin{figure}[ht]
        \begin{minipage}[b]{0.325\linewidth}
            \centering
            \subfloat[\footnotesize{Case 1}\label{fig:same_cap}]{%
i                \includegraphics[width=1\linewidth]{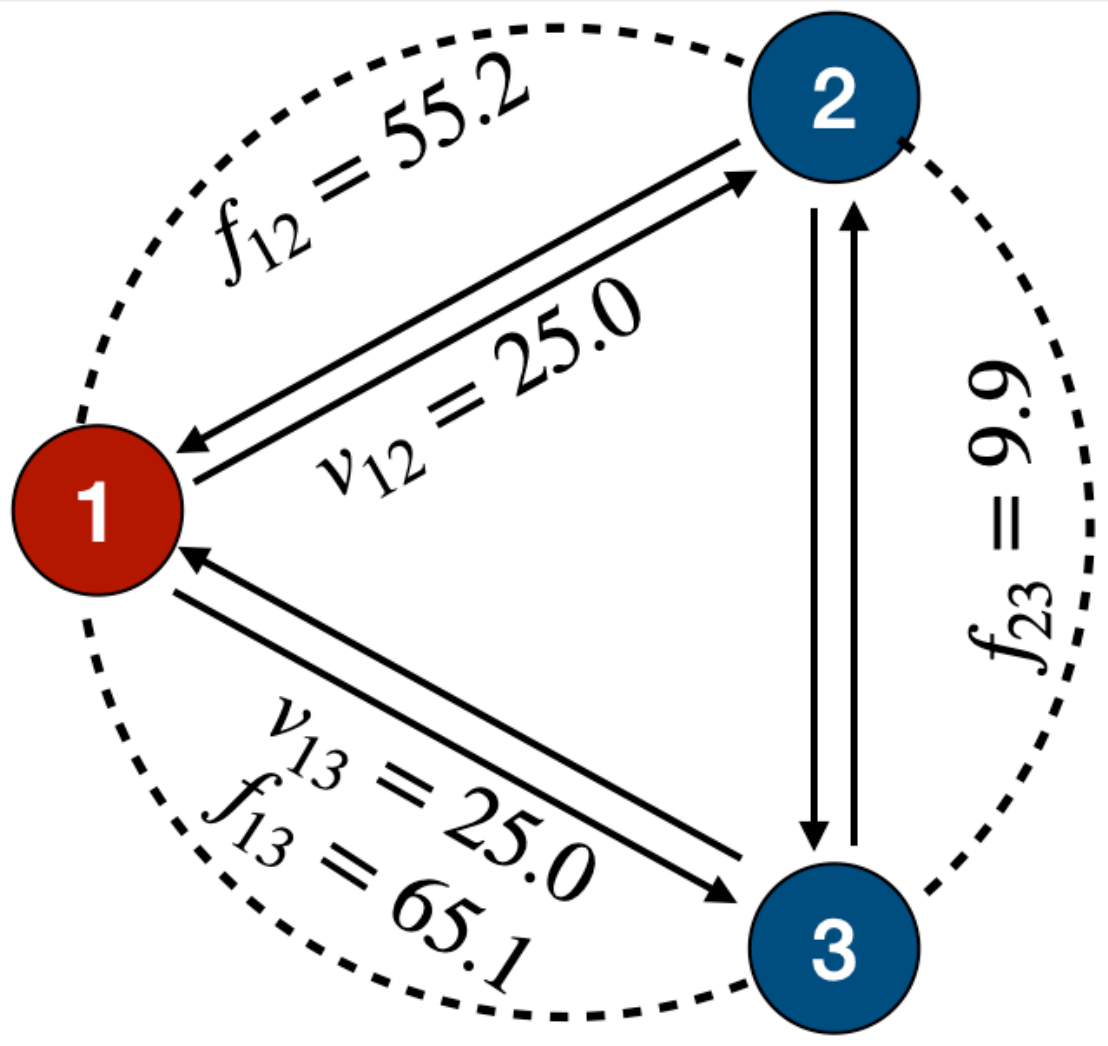}}
        \end{minipage}
        \begin{minipage}[b]{0.325\linewidth}
            \centering
            \subfloat[\footnotesize{Case 2}\label{fig:road_congestion_scen}]{%
                \includegraphics[width=1\linewidth]{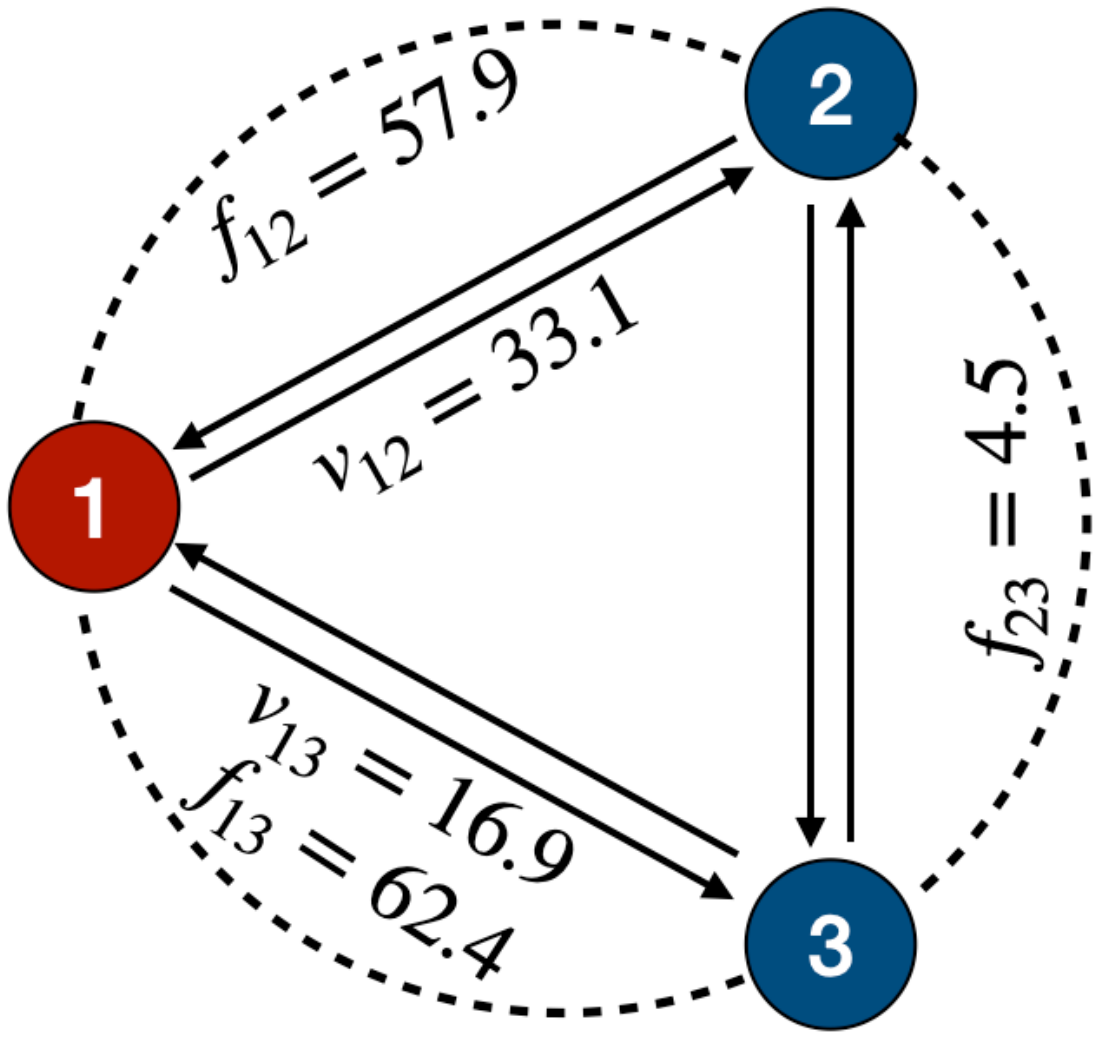}}
        \end{minipage}
        \begin{minipage}[b]{0.325\linewidth}
            \centering
            \subfloat[\footnotesize{Case 3}\label{fig:power_congestion_scen}]{%
                \includegraphics[width=1\linewidth]{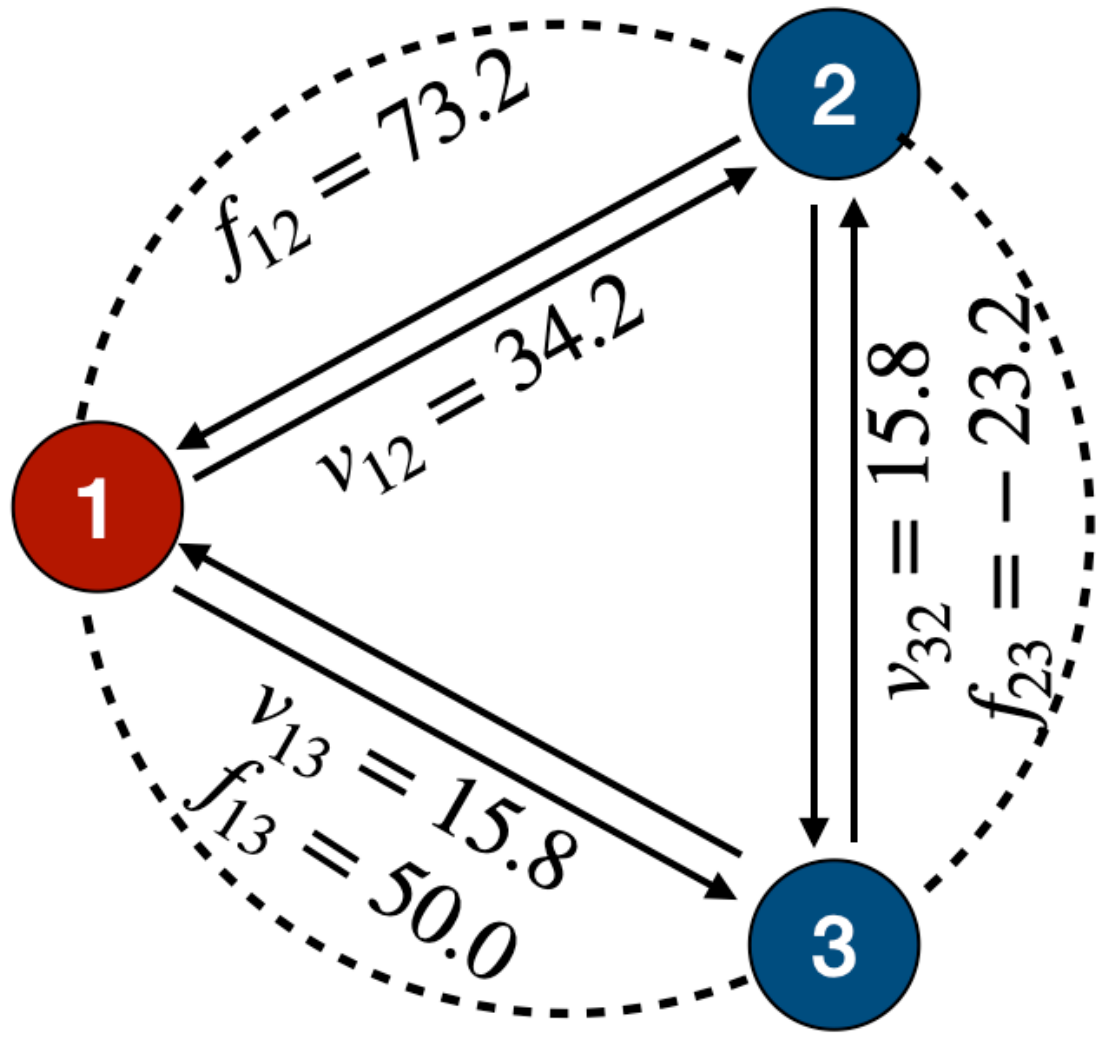}}
        \end{minipage}
        \caption{Impacts of link capacity on link flow.}
        \label{fig:capacity_impact}
    \end{figure}

	Comparing between Fig. \ref{fig:same_cap} and Fig. \ref{fig:road_congestion_scen}, transportation congestion on link 1-3 shifts 8.1 units of travel demand from link 1-3 to link 1-2. The increasing of charging demand on node 2 leads to redistribution of power flow on each transmission line to preserve power flow conservation and power physics laws. Similarly, comparing between Fig. \ref{fig:same_cap} and Fig. \ref{fig:power_congestion_scen},  constraining the transmission capacity of link 1-3 results an increase of the power flow on links 1-2 and 2-3. In addition, the reduced transmission capacity on link 1-3 will increase the energy prices at node 3 more than the prices at node 2 (see Section \ref{sec:equi_price}), which will discourage EVs from traveling to node 3. Notice that few vehicles are charging at node 3 in case 3 (see Fig. \ref{fig:power_congestion_scen}), and the traffic flow from node 1 to node 2 splits into two paths, $1-2$ and $1-3-2$, to avoid traffic congestion. \revi{Again, these results also illustrate the interdependence between transportation and power systems. Without properly considering the complicated interaction and feedback effects, the analyses results may be biased.}

	
    



    \subsubsection{Sensitivity Analysis of Load and EV Penetration}
    
    \revi{We have investigated the sensitivity of existing power load and EV penetration on the system. The results for both case 1 and 2 showed similar patterns in energy price and link traffic flow except that the traffic flow splits unevenly between the links in case 2. Therefore, we have presented the results only for case 2 and 3 here to maintain conciseness.}   
    

     
    \revi{Firstly, we change the existing load demand from 20 to 50 in the increments of 5. In case 2, the energy price increases linearly in all three nodes (see Fig. \ref{fig:case2_energy_price_change_load}). Since we have the same price for all the nodes and all the scenarios in each increment of load change, both nodes 2 and 3 are equally desirable for charging in terms of charging cost. However, the road capacity of link 1-3 is less than that of link 1-2, and the traffic flow splits unevenly between these two links (see Fig. \ref{fig:case2_traffic_change_load}).
    
    We reach more intriguing results for case 3; the energy prices start increasing for the three nodes until the load reaches 50 (Fig. \ref{fig:case3_price_change_load}). After that, the problem becomes infeasible, which is sensible because node 3 can not receive its required energy from node 1 and node 2 due to transmission capacity limitation (\ref{cons:iso_capa}) and phase angle constraints (\ref{cons:iso_phas}). Furthermore, node 3 has a higher charging price as the locational load increases than node 2. Therefore, the traffic flow to node 2 is typically more than the flow to node 3. While most of the drivers use link 1-2 to reach node 2, some drivers choose links 1-3 and 3-2 to avoid congestion on link 1-2 (see Fig. \ref{fig:case3_traffic_change_load}).}
        \revi{
    	\begin{figure}[ht]
        \begin{minipage}[b]{0.5\linewidth}
            \centering
            \subfloat[\footnotesize{Case 2}\label{fig:case2_energy_price_change_load}]{%
                \includegraphics[width=0.8\linewidth]{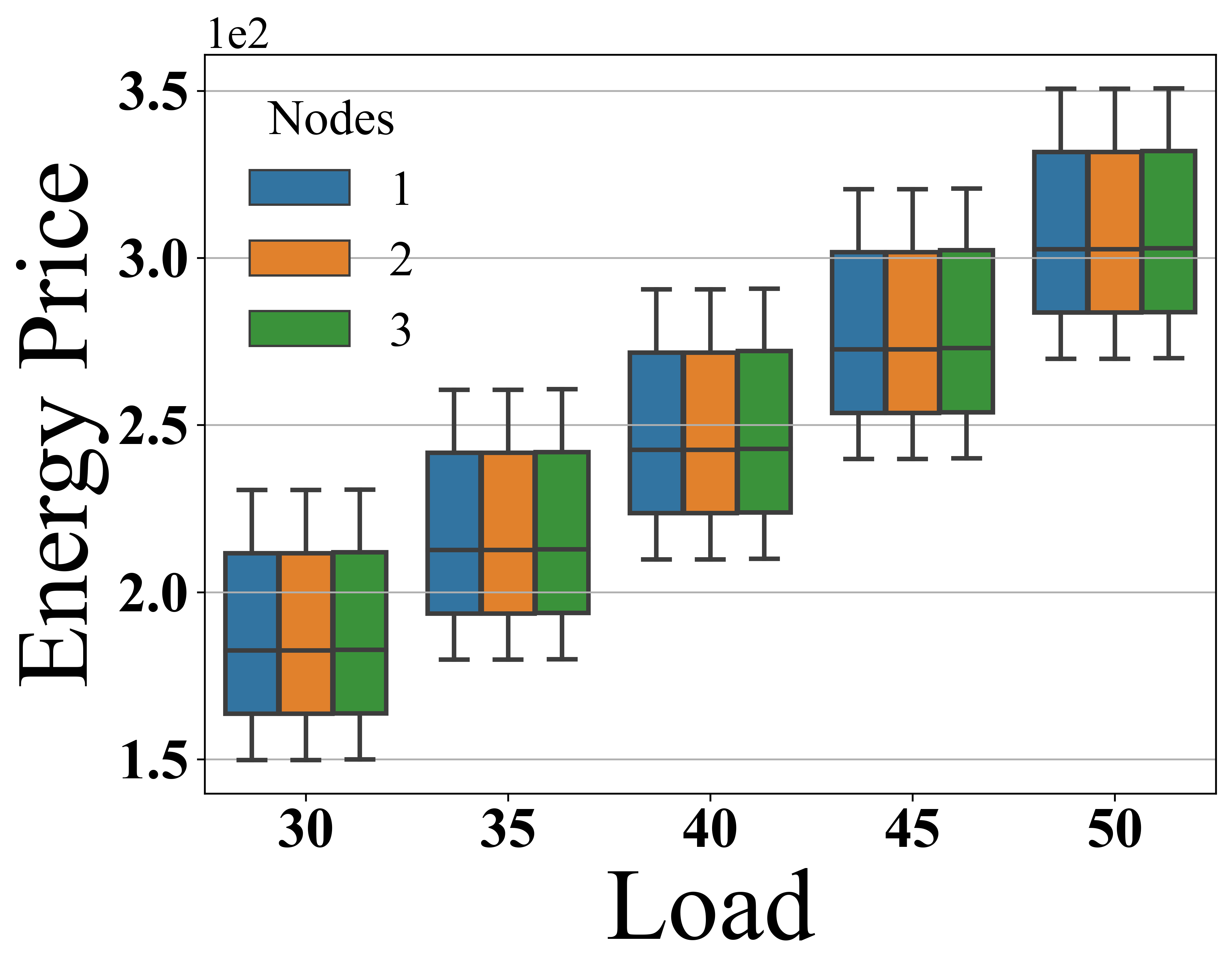}}
        \end{minipage}
        \begin{minipage}[b]{0.5\linewidth}
            \centering
            \subfloat[\footnotesize{Case 3}\label{fig:case3_price_change_load}]{%
                \includegraphics[width=0.8\linewidth]{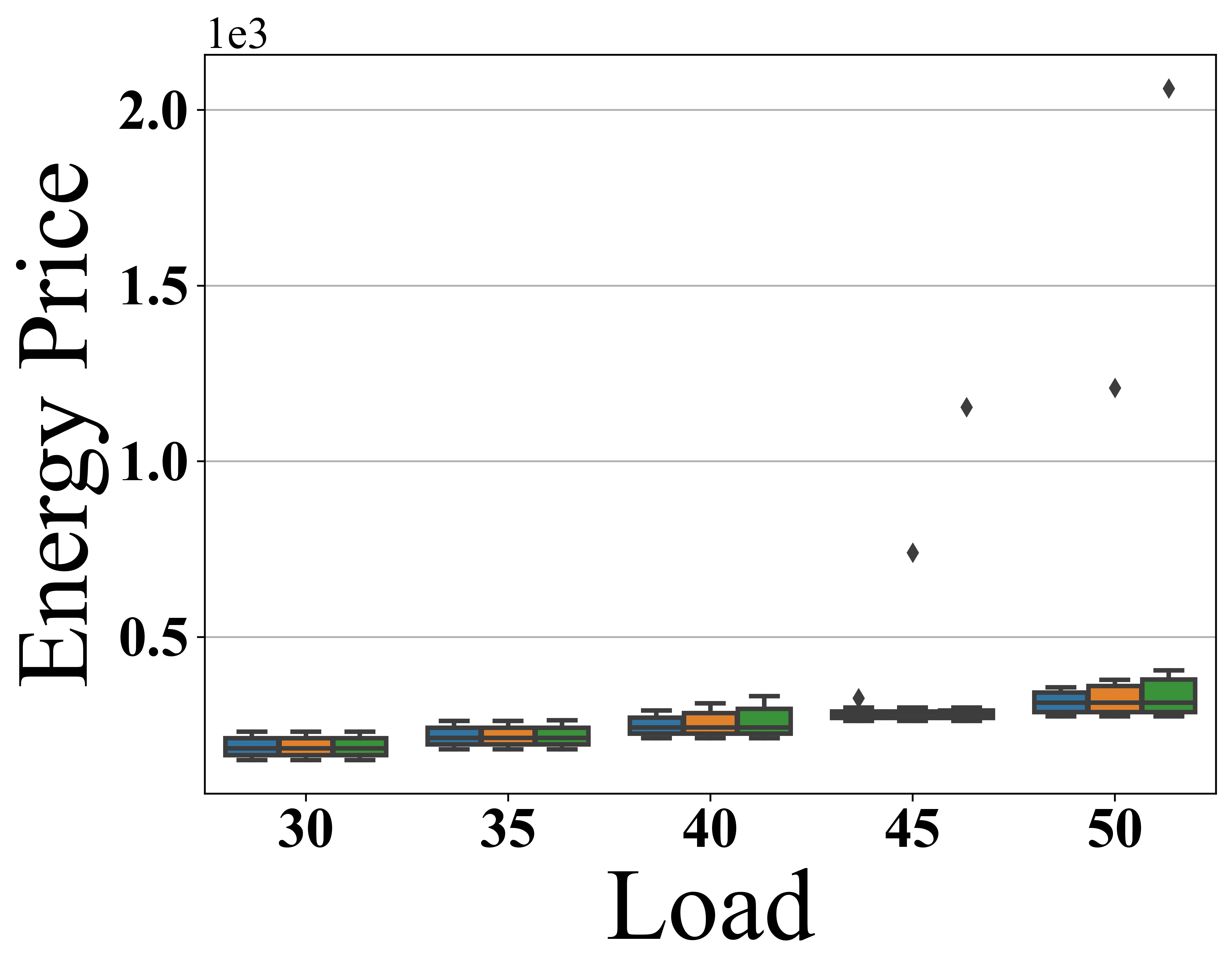}}
        \end{minipage}
        
        \begin{minipage}[b]{0.5\linewidth}
            \centering
            \subfloat[\footnotesize{Case 2}\label{fig:case2_price_change_EV}]{%
                \includegraphics[width=0.8\linewidth]{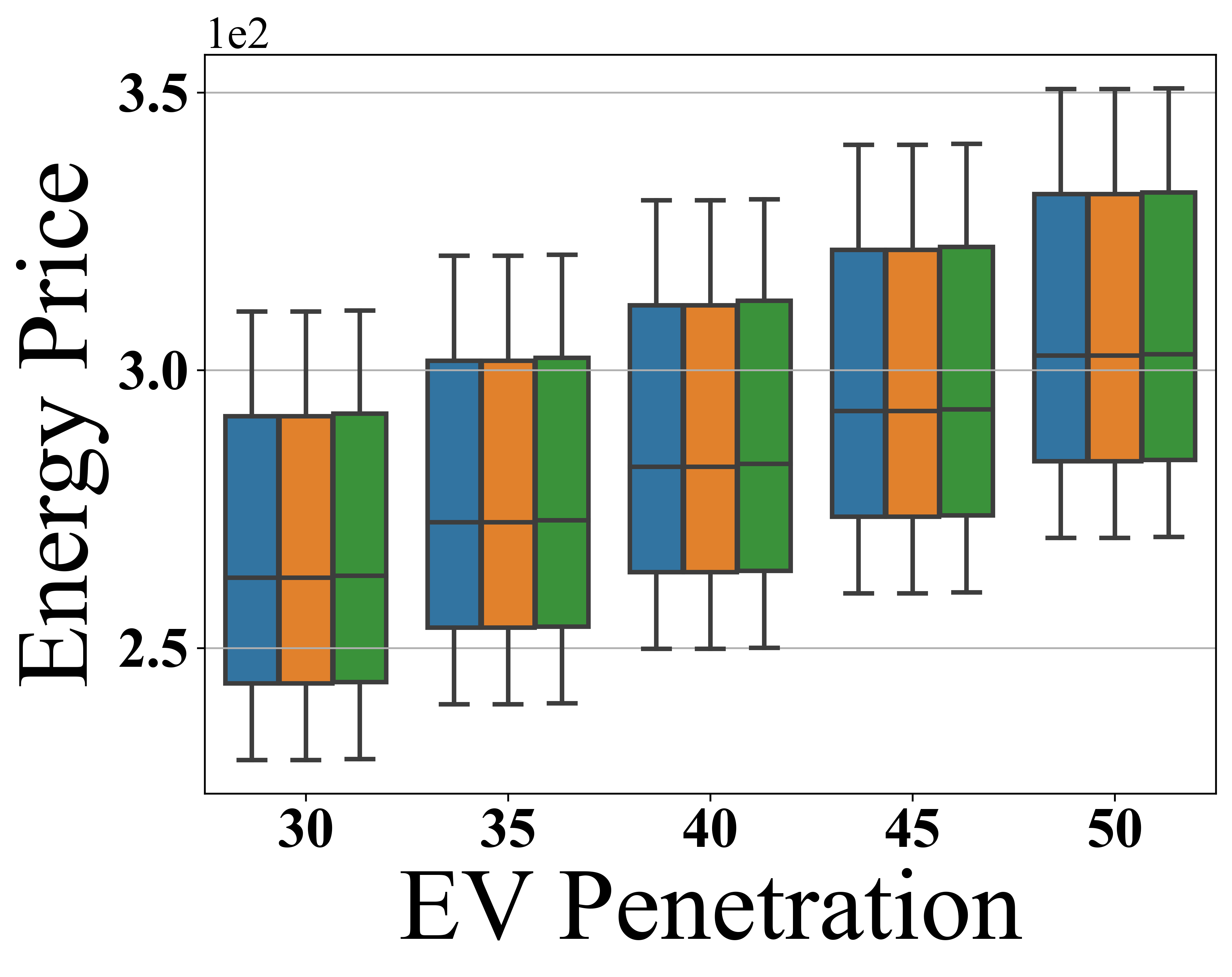}}
        \end{minipage}
        \begin{minipage}[b]{0.5\linewidth}
            \centering
            \subfloat[\footnotesize{Case 3}\label{fig:case3_price_change_EV}]{%
                \includegraphics[width=0.8\linewidth]{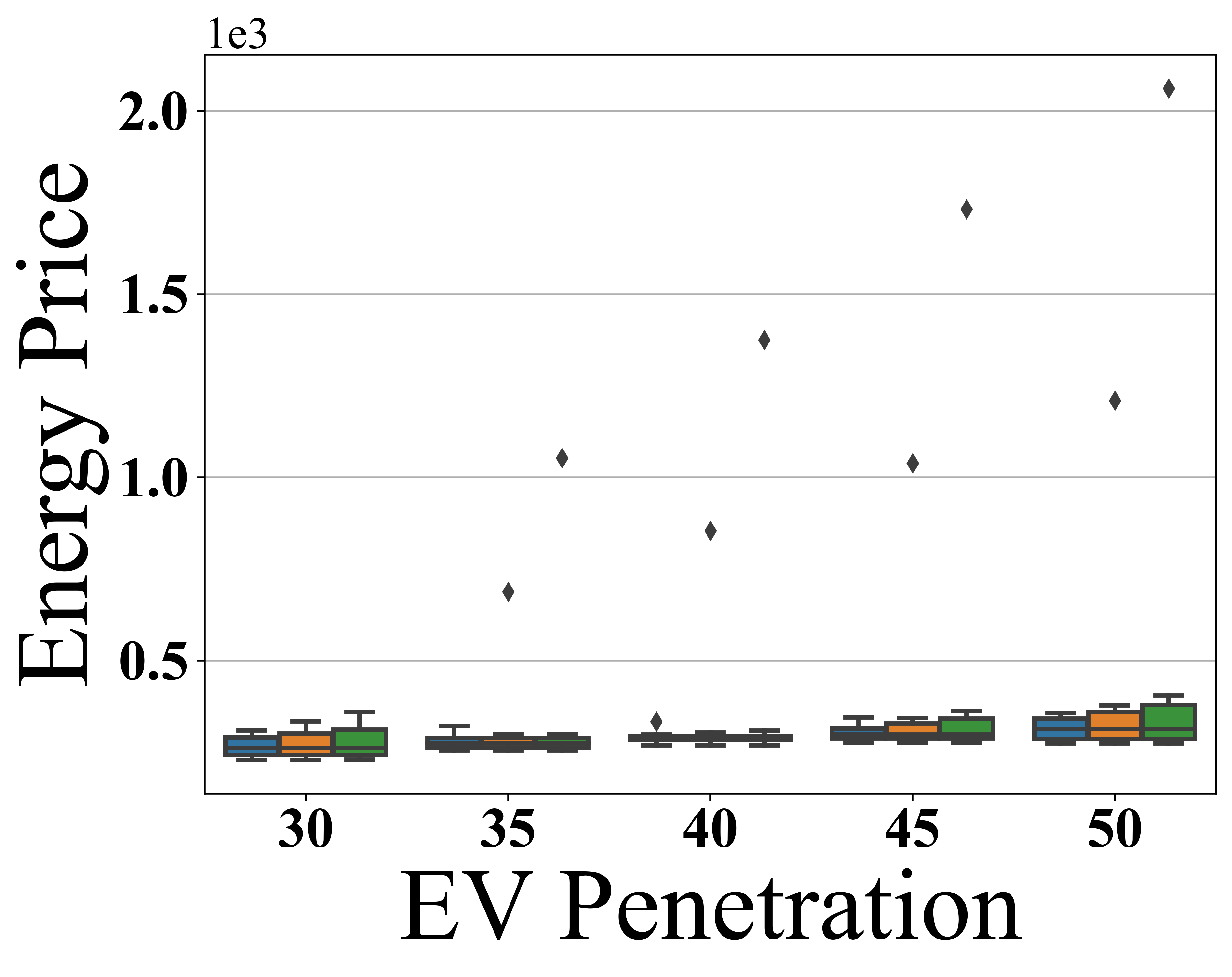}}
        \end{minipage}

        \caption{Impacts of existing load demand and EV penetration on Energy Price.}
        \label{fig:Load_change_impact_case2}
    \end{figure}
    }
    \revi{
    	\begin{figure}[ht]
        \begin{minipage}[b]{0.5\linewidth}
            \centering
            \subfloat[\footnotesize{Case 2}\label{fig:case2_traffic_change_load}]{%
                \includegraphics[width=0.8\linewidth]{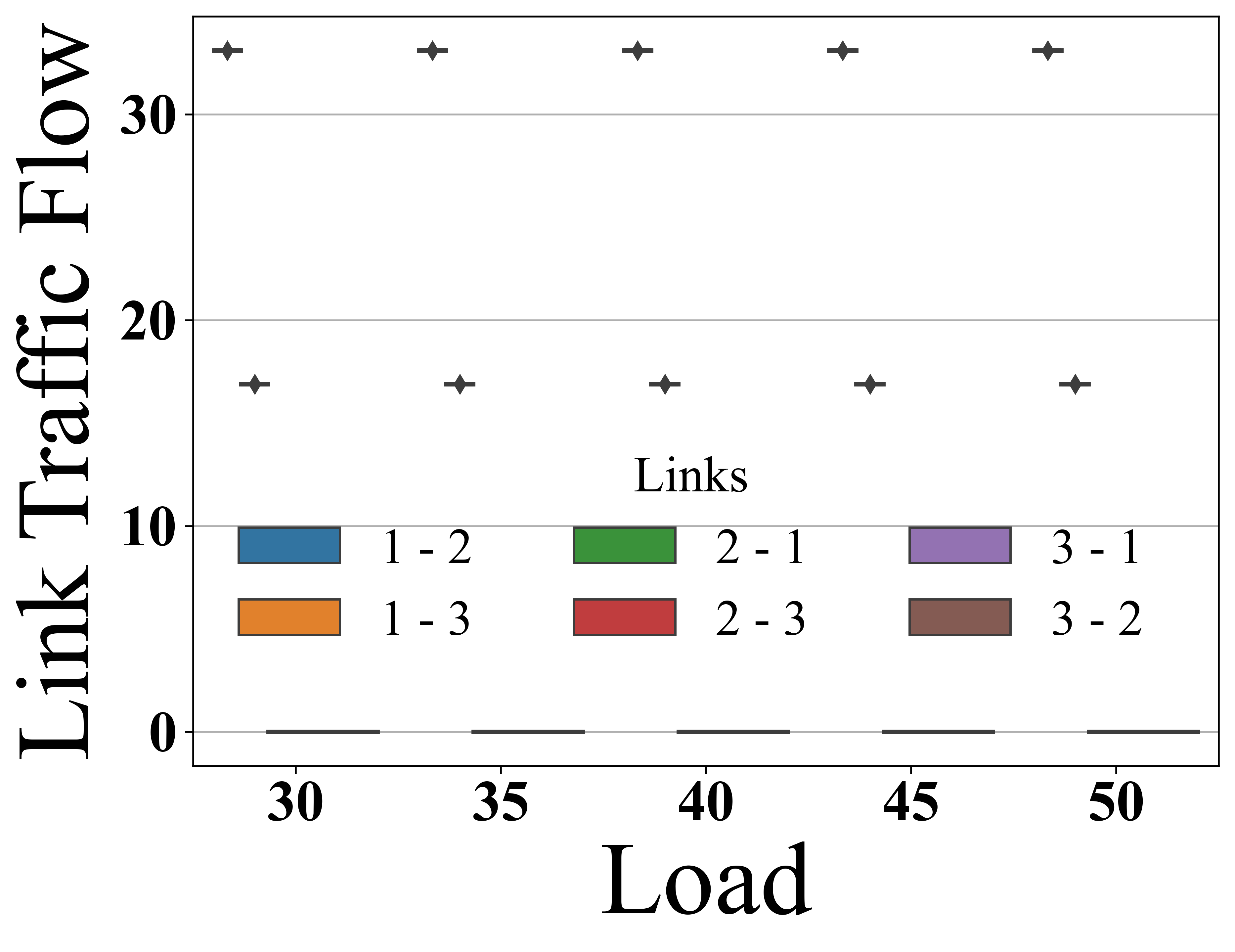}}
        \end{minipage}
        \begin{minipage}[b]{0.5\linewidth}
            \centering
            \subfloat[\footnotesize{Case 3}\label{fig:case3_traffic_change_load}]{%
                \includegraphics[width=0.8\linewidth]{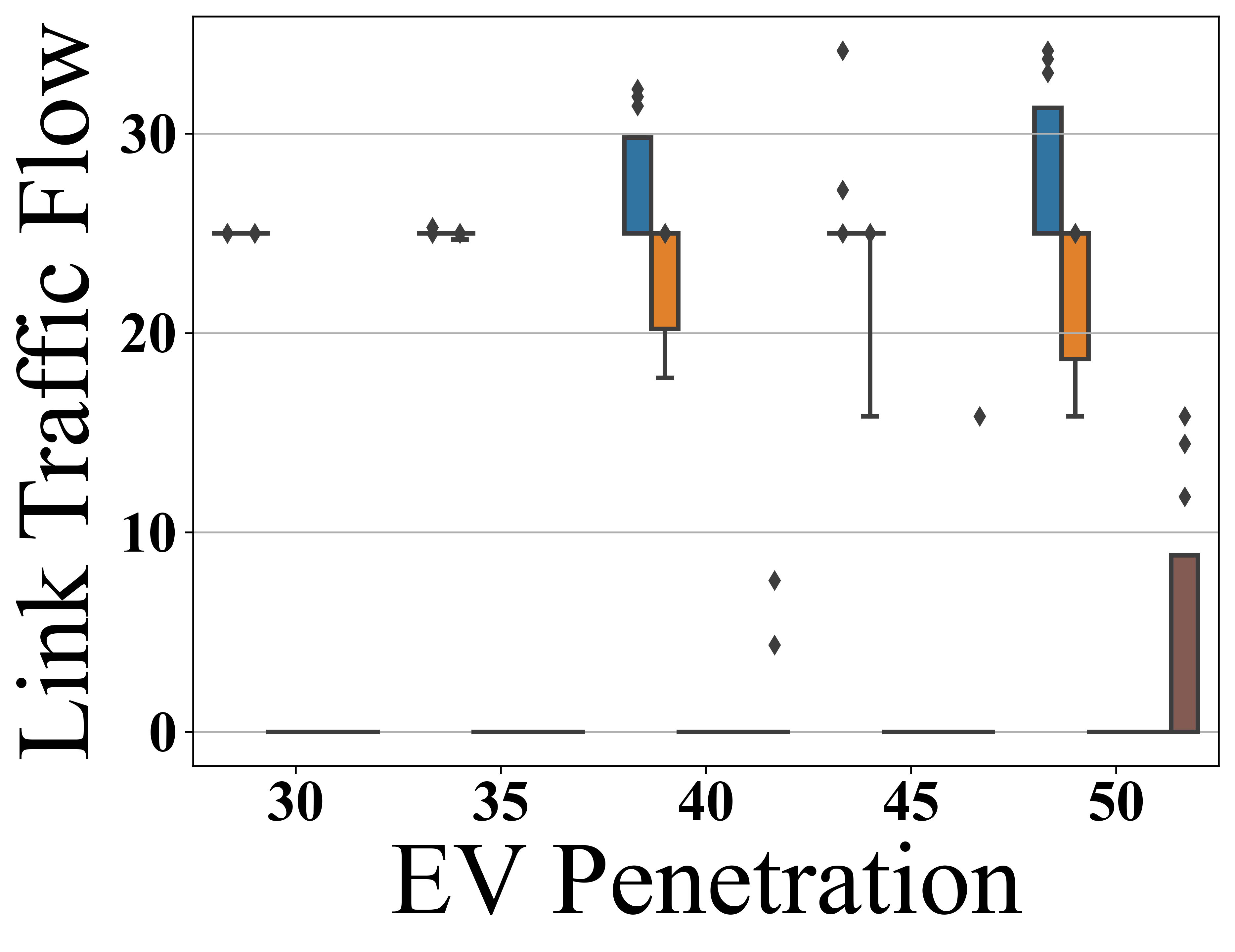}}
        \end{minipage}
        
        \begin{minipage}[b]{0.5\linewidth}
            \centering
            \subfloat[\footnotesize{Case 2}\label{fig:case2_traffic_change_EV}]{%
                \includegraphics[width=0.8\linewidth]{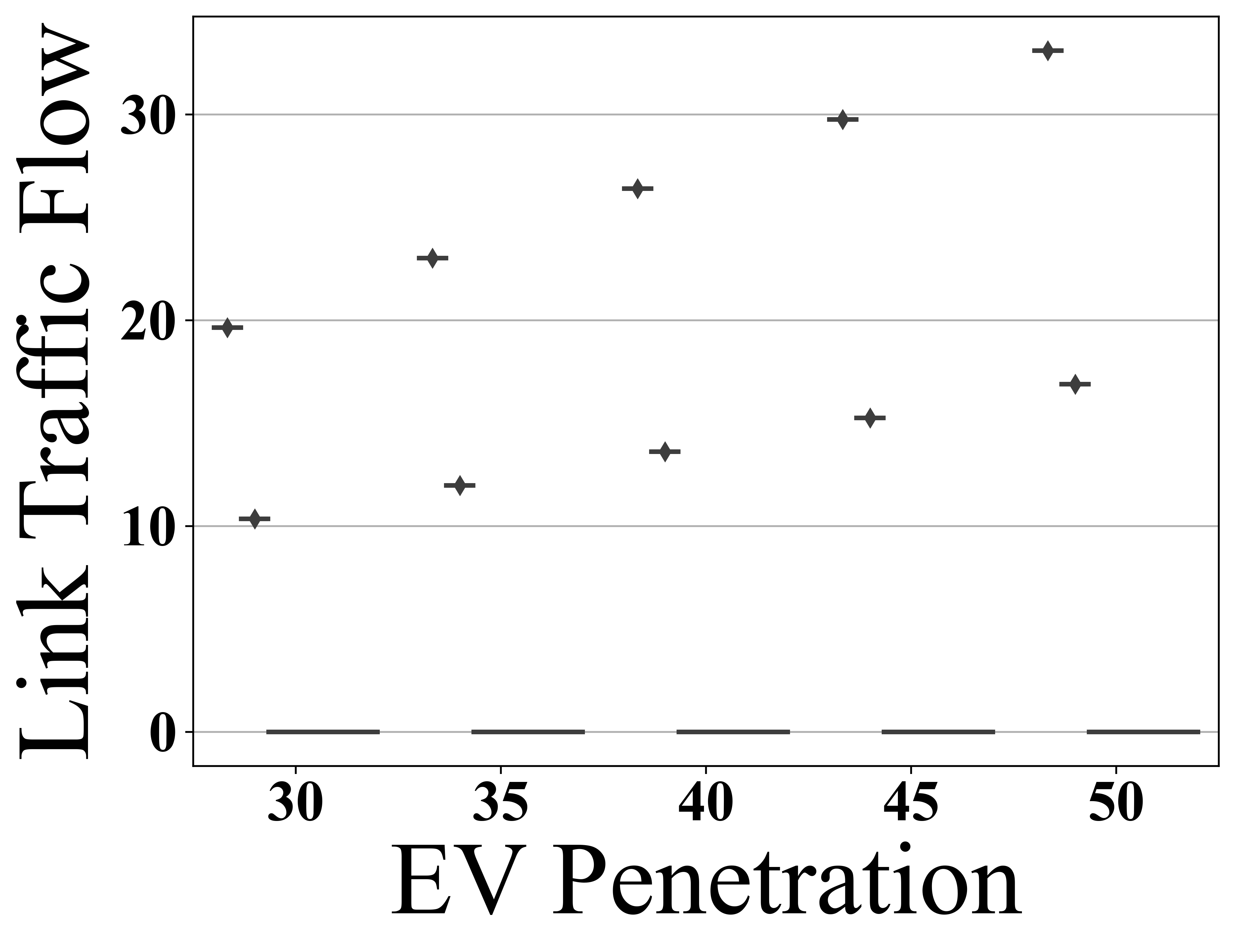}}
        \end{minipage}
        \begin{minipage}[b]{0.5\linewidth}
            \centering
            \subfloat[\footnotesize{Case 3}\label{fig:case3_traffic_change_EV}]{%
                \includegraphics[width=0.8\linewidth]{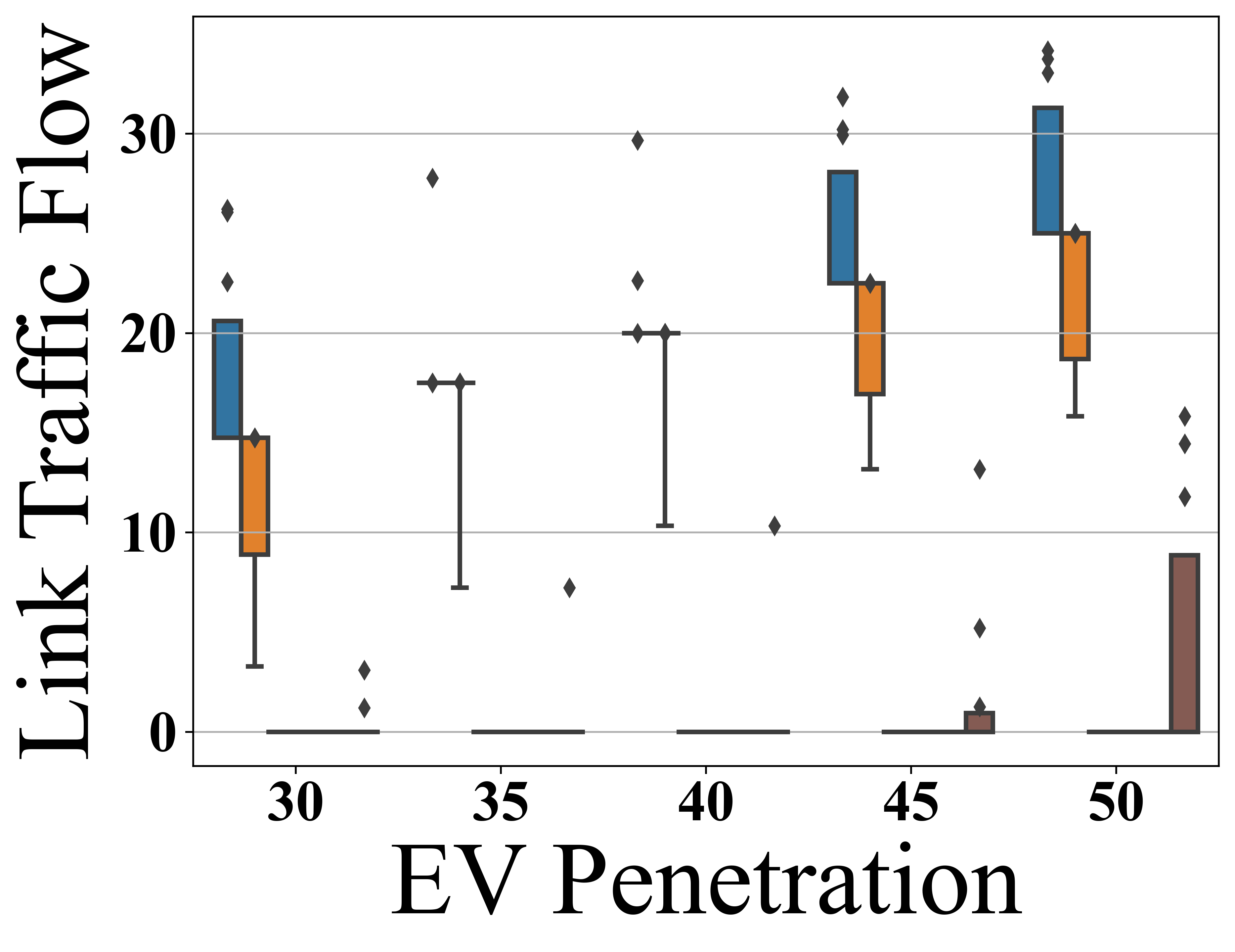}}
        \end{minipage}

        \caption{Impacts of existing load demand and EV penetration on link traffic flow.}
        \label{fig:Link_flow_analysis}
    \end{figure}
    }
    
    \revi{We also change the initial EV traffic flow from 20 to 50 with increments of 5 with the load demand to be 50 in all nodes. 
    As the EV penetration increases, the traffic flow on both links 1-2 and 1-3 also increases.  With limited road capacity on link 1-3 in case 2, however, we see more traffic on link 1-2 compared to link 1-3 (see Fig. \ref{fig:case2_price_change_EV}). On the other hand, the energy prices increase linearly for all nodes, because the EV penetration acts as load demand on power system (see Fig. \ref{fig:case2_traffic_change_EV}). For case 3, the limitation on transmission line 1-3 changed the linear increase in energy prices (Fig. \ref{fig:case3_price_change_EV}), and similar to the case 3 for the load change, the optimization problem becomes infeasible when EV penetration exceeds 50. The difference on the locational energy prices has also caused different route selection of EVs, changing the dynamic of traffic flow on the links (see Fig. \ref{fig:case3_traffic_change_EV}).
    
    The results found for sensitivity analysis of both EV penetration and existing load demand show the necessity to model the interdependence of transportation and power systems, where limitations on one system may cause drastic impacts on the other.} 

    \vspace{-0.5cm}
    \subsection{Sioux Falls Road Network and IEEE 39-bus Test System}
\begin{figure}[ht]
    	\begin{minipage}[b]{0.55\linewidth}
            \centering
            \includegraphics[angle=90,origin=c,width=1\linewidth]{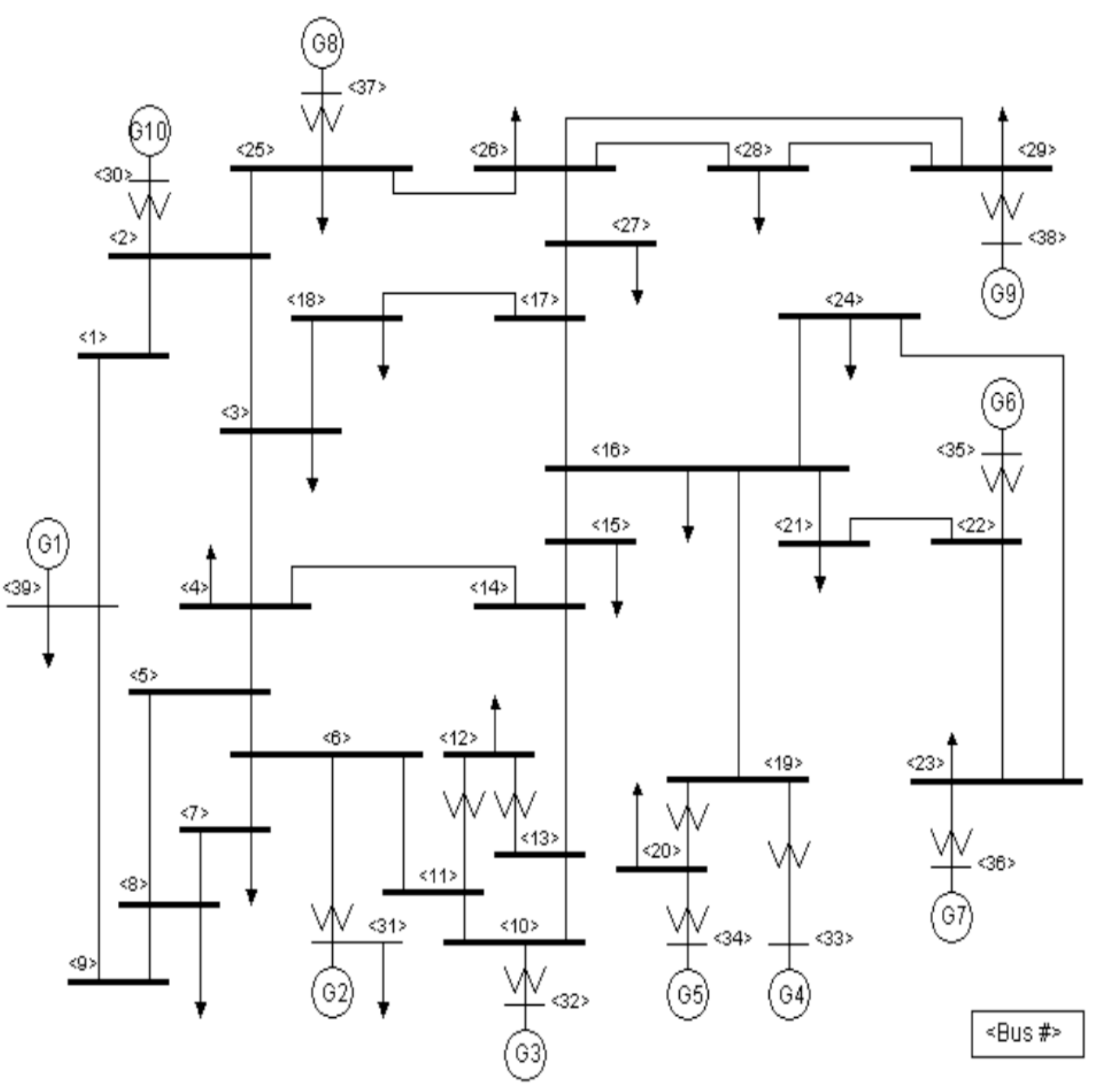}
            \caption{IEEE 39-bust test systems.}
            \label{fig:39_bus}
        \end{minipage}
        \begin{minipage}[b]{0.44\linewidth}
            \centering
            \includegraphics[width=1\linewidth]{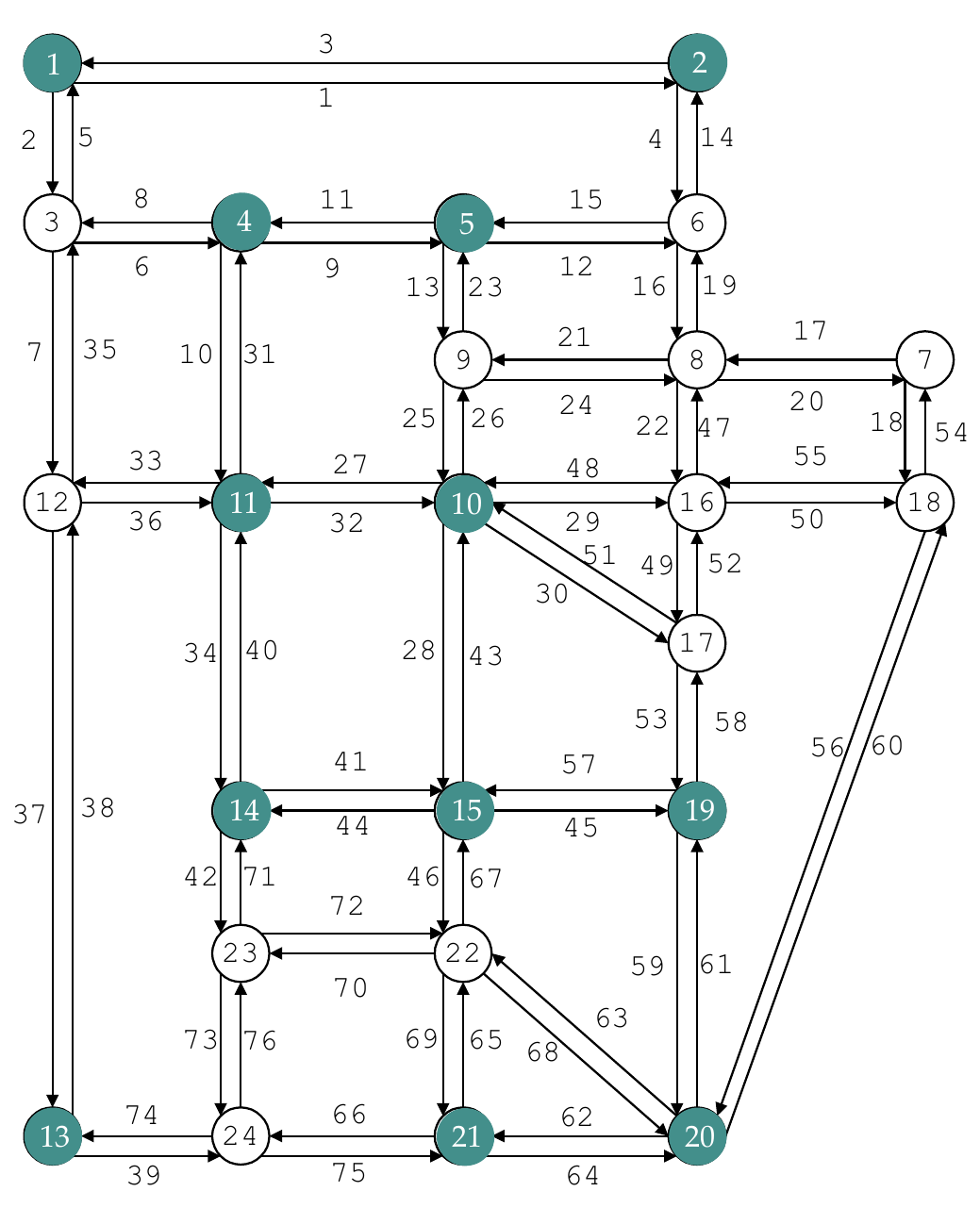}
            \caption{Sioux falls transportation network.}
            \label{fig:siou_fall}
        \end{minipage}
	\end{figure}

\begin{figure}[!t]
		\centering
		\includegraphics[width=1\linewidth]{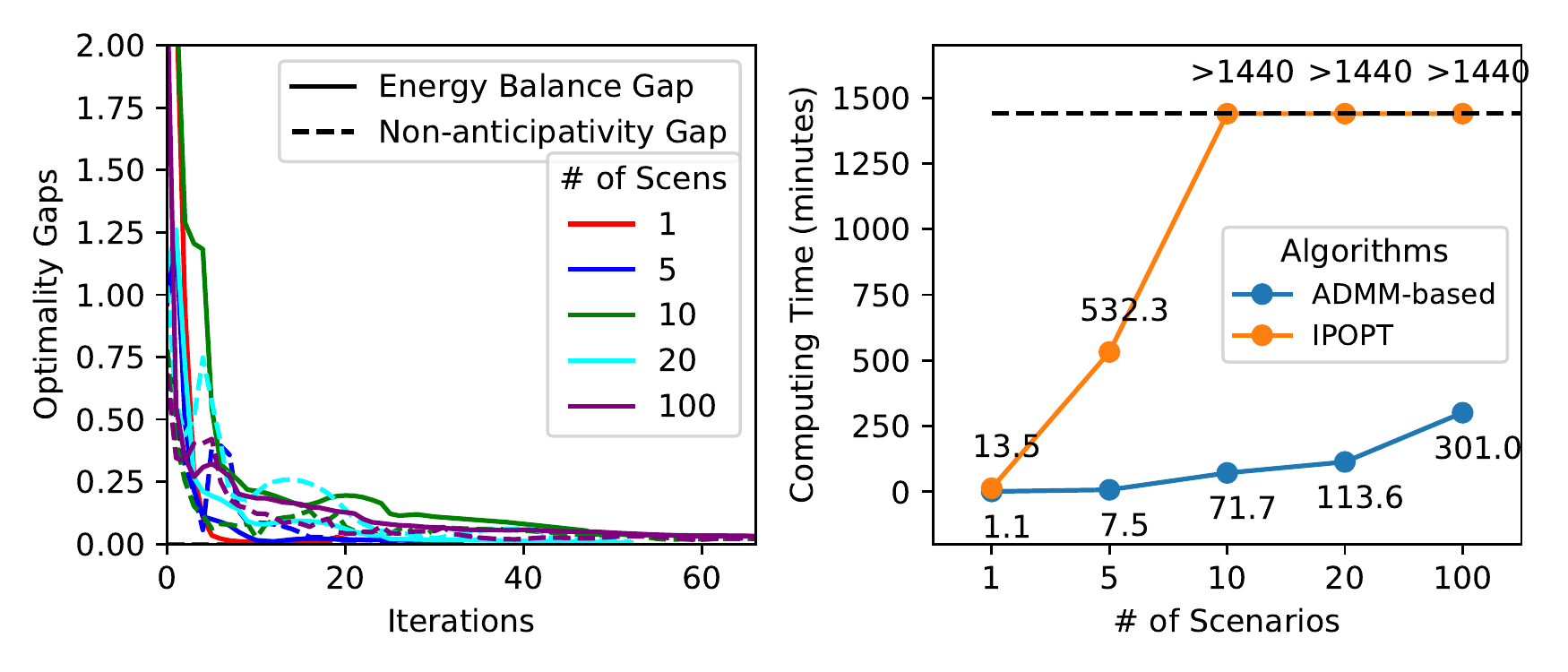}
		\caption{Convergence patterns and computing time.}
		\label{fig:convergence}
	\end{figure}

    \rev{In addition to insights generated by small test cases above, this section aims to demonstrate the convergence properties of our algorithms on larger networks. We use Sioux Falls road network\footnote{Data: \url{https://github.com/bstabler/TransportationNetworks/}} (see Fig. \ref{fig:siou_fall}) and IEEE 39-bus test system\footnote{Data: \url{https://matpower.org/docs/ref/matpower5.0/case39.html}} (see Fig. \ref{fig:39_bus}), two test systems widely used in transportation and power system literature, respectively.} The correspondence between the node indexes in transportation and power systems is shown in Table \ref{tbl:node_corr}. To avoid infeasibility due to increased power demand from charging. we scale down the travel demand and road capacity to be 1\% of the original values. \revi{Sioux Falls network is at city level, while IEEE 39-bus system is a regional-level transmission network. To better reflect the inter-city travel, we scale up the free flow link travel time in Sioux Falls network to 10 times of its original value. A similar approach was adopted by \cite{he2013optimal, guo2016infrastructure}. The computational results presented here  are for illustration purposes only. Future research may need to collect transportation and power network data from the same geographical area to gain practical insights.} \rev{The candidate EVs charging locations and PV investment locations are marked in green in Fig. \ref{fig:siou_fall}}. The charging load accounts for 18.5\% of the total load. The uncertain renewable generation factors $\boldsymbol{\xi}$ are randomly sampled from uniform distribution $[0.5, 1.5]$. The optimality gap is set to be 1\%.
	\begin{table}[ht]
 	\caption{Node Correspondence between Systems}
	\label{tbl:node_corr}
	\centering
	\resizebox{0.8\columnwidth}{!}{
		\begin{tabular}{ccccccccccccc}
            \hline
            System& \multicolumn{12}{c}{Node Index}\\
			\hline
			Transportation	& 1	&2	&4	&5	&10	&11	&13	&14	&15	&19	&20	&21\\
			Power	&1	&4	&6	&11&	13&	16&	19&	2&	23&	25&	27&	32\\
			\hline
		\end{tabular}
	}
	\end{table}%
	The convergence patterns and computing time are shown in Fig. \ref{fig:convergence}. The algorithm we developed converges reliability within \revi{100} iterations for up to \revi{100} scenarios. The computing time is almost linearly increasing (from \revi{1.1} minutes to \revi{301.0} minutes) with the scenario number for our algorithm, but with a higher increasing rate when solving the whole problem using IPOPT. IPOPT cannot solve for 10 or more scenarios in 24 hours. \revi{We note that our algorithm has the potential for more scenarios due to the solution strategies of scenario decomposition and parallel computing. But one may need to increase the number of CPU cores to full materialized the benefits brought by parallel computing. In addition, since our main computation strategies decompose the transportation and power systems, the algorithm could be further optimized by using more advanced algorithms for transportation network assignment and solving power flow equations. These are beyond the scope of this paper and will be left for the future.}

    \section{Conclusion} \label{sec:conclusion}
	We study the interdependence between transportation and power systems considering decentralized renewable generation and electric vehicles. We propose a stochastic multi-agent optimization framework to formulate the complex interactions between  stakeholders in transportation and power systems, including EV/conventional vehicle drivers, \revi{renewable}/conventional generators, and independent system operators, with locational electricity and charging prices endogenously determined by markets. We develop efficient computational algorithms to cope with the curse of dimensionalities based on exact convex reformulation and decomposition-based approaches. We prove the existence and uniqueness of the equilibrium. Numerical experiments show that our algorithms outperform the existing commercial solvers, and provide insights on the close coupling between transportation and power systems. 

This work can be extended in several directions. From a modeling perspective, we adopted classic transportation and power flow models. Investigating other alternatives, such as traffic equilibrium models with intermediate stops, may generate more accurate descriptions of traffic flow. Second, our current model could be extended to consider other stochastic factors in addition to renewable generation, such as load, travel, and charging demand uncertainties. But more sophisticated models and high dimensional stochasticity may be more computational demanding to solve. \revi{Third, we assume a perfectly competitive market for renewable energy investment. It is valuable to investigate different market settings in order to gain a better understanding about the system interactions.} From a computation perspective, adopting more advanced solution algorithms for transportation and power systems will further speed up the solution for subproblems. In addition, effectively grouping uncertain scenarios will accelerate the convergence rate. \revii{From an application perspective, this model can be served as a foundation to capture the interdependence between transportation and power systems, upon which further system planning, operation, restoration, and incentive design strategies can be modeled and investigated to improve efficiency, reliability and resilience of both systems.}
	
	
	%

	\appendices
	
	\section{}
\label{app:proofs}	
	\IEEEproof (Theorem \ref{thm:refo_perf}) The Lagrangian of (\ref{mod:eq_perf}) can be written as (\ref{lag:eq_perf}) after relaxing equilibrium constraints (\ref{eq:equi}).
	
	\resizebox{1\columnwidth}{!}{%
	\begin{minipage}{\columnwidth}
		\begin{align}
		\mathcal{L} &= \sum_{i \in \mathcal{I}^S} C_i^{S,I}(u_i^S) +  \mathbb{E}_{\boldsymbol{\xi}} \bigg[ \sum_{i \in \mathcal{I}^S}C_i^{S,O}(g_{i,\boldsymbol{\xi}}^S)) +  \sum_{i \in \mathcal{I}^C} C^C_i(g_{i,\boldsymbol{\xi}}^C) \nonumber \\
		&+ \frac{\beta_1}{\beta_2}\sum_{a \in \mathcal{A}} \int_{0}^{v_{a,\boldsymbol{\xi}}} tt_a(u) \mathrm{d}u + \frac{1}{\beta_2}\sum_{r \in \mathcal{R}}\sum_{s \in \mathcal{S}} q_{rs,\boldsymbol{\xi}}\left(\ln q_{rs,\boldsymbol{\xi}} - 1 - \beta_{0,s}\right) \nonumber\\
		& + \sum_{i \in \mathcal{I}^S \cup \mathcal{I}^C }\tilde{\rho}_{i,\boldsymbol{\xi}}(d_{i,\boldsymbol{\xi}} - g_{i,\boldsymbol{\xi}}^S - g_{i,\boldsymbol{\xi}}^C) + \sum_{s \in \mathcal{S}}\tilde{\lambda}_{s,\boldsymbol{\xi}}\Big(\sum_{r \in \mathcal{R}}e_{rs}q_{rs,\boldsymbol{\xi}} - p_{i(s),\boldsymbol{\xi}}\Big)\bigg]\nonumber
		\end{align}
		\end{minipage}}
		\resizebox{1\columnwidth}{!}{%
	\begin{minipage}{\columnwidth}
		\begin{align}
		&= \sum_{i \in \mathcal{I}^S} C_i^{S,I}(u_i^S)+ \mathbb{E}_{\boldsymbol{\xi}} \bigg[\sum_{i \in \mathcal{I}^S}C_i^{S,O}(g_{i,\boldsymbol{\xi}}^S) - \tilde{\rho}_{i,\boldsymbol{\xi}} g_{i,\boldsymbol{\xi}}^S \nonumber\\
		&+ \sum_{i \in \mathcal{I}^C}C_i^{C}(g_{i,\boldsymbol{\xi}}^C) - \tilde{\rho}_{i,\boldsymbol{\xi}} g_{i,\boldsymbol{\xi}}^C + \frac{\beta_1}{\beta_2}\sum_{a \in \mathcal{A}} \int_{0}^{v_{a,\boldsymbol{\xi}}} tt_a(u) \mathrm{d}u \nonumber\\
		&+ \frac{1}{\beta_2}\sum_{r \in \mathcal{R}}\sum_{s \in \mathcal{S}} q_{rs,\boldsymbol{\xi}}\left(\ln q_{rs,\boldsymbol{\xi}} - 1 +\beta_2\tilde{\lambda}_{s,\boldsymbol{\xi}}e_{rs} - \beta_{0,s}\right) \nonumber\\
		&+\sum_{i \in \mathcal{I}^S \cup \mathcal{I}^C }\tilde{\rho}_{i,\boldsymbol{\xi}}d_{i,\boldsymbol{\xi}} + \sum_{s \in \mathcal{S}}\tilde{\lambda}_{s,\boldsymbol{\xi}}(- p_{i(s),\boldsymbol{\xi}})\bigg]
		\label{lag:eq_perf}
		\end{align}
		\end{minipage}}
	Define feasible set $\mathcal{X} \doteq \{\boldsymbol{u,g,p,d,x,q}\geq \boldsymbol{0} |(\ref{cons:sp_capa})-(\ref{cons:sp_budg}), (\ref{cons:cg_capa})-(\ref{cons:cg_lower_bound}), (\ref{cons:iso_phas})-(\ref{cons:iso_capa}), (\ref{cons:cda_v_x})-(\ref{cons:cda_q_d})\}$. 
	Models (\ref{mod:sp}), (\ref{mod:cp}), (\ref{mod:iso}), (\ref{mod:cda}), and (\ref{mod:eq_perf}) are all convex optimization and satisfy linearity constraint qualifications. So strong duality holds. (\ref{mod:eq_perf}) is equivalent to (\ref{mod:dual}) due to strong duality.
	\begin{equation}\label{mod:dual}
	     \min_{(\boldsymbol{u,g,p,d,x,q}) \in \mathcal{X}}\max_{\boldsymbol{\tilde{\lambda}, \tilde{\rho}}} \mathcal{L} = \max_{\boldsymbol{\tilde{\lambda}, \tilde{\rho}}} \min_{(\boldsymbol{u,g,p,d,x,q}) \in \mathcal{X}} \mathcal{L}
	\end{equation}
	$\min_{(\boldsymbol{u,g,p,d,x,q}) \in \mathcal{X}} \mathcal{L}$ is equivalent with (\ref{mod:sp}, \ref{mod:cp}, \ref{mod:iso}, \ref{mod:cda}) for any given $\boldsymbol{\tilde{\lambda}, \tilde{\rho}}$. In addition, (\ref{eq:equi}) holds for the optimal $\boldsymbol{\tilde{\lambda}, \tilde{\rho}}$. \IEEEQED

    \IEEEproof (Corollary \ref{cor:exis_uniq_syst_equi}) If $tt_a(\cdot)$, $C_i^{S,I}(\cdot)$, $C_i^{S,O}(\cdot)$, and $C_i^{C}(\cdot)$  are strictly convex functions, model \ref{mod:eq_perf} is strict convex optimization problem, which has a unique optimal solution. Following Theorem \ref{thm:refo_perf}, the system equilibrium therefore exists and is unique. \IEEEQED


	\ifCLASSOPTIONcaptionsoff
	\newpage
	\fi

	
	
	
	
	\bibliographystyle{IEEEtran}
	\bibliography{IEEEabrv,./my_library.bib}
	%
	%
	%
	%
	%
\end{document}